\documentclass[11pt]{article}
\usepackage{amsfonts,mathrsfs,amssymb,amsthm,mathptm}
\usepackage{amsmath,amscd}
\usepackage{mathptm,pslatex}
\usepackage{color}
\usepackage[dvips]{graphicx}
\usepackage[all]{xy}
\usepackage{graphicx}
\usepackage{fancyhdr}
\usepackage{url}

\begin{document}
\newtheorem{defi}{Definition}[section]
\newtheorem{exam}[defi]{Example}
\newtheorem{prop}[defi]{Proposition}
\newtheorem{theorem}[defi]{Theorem}
\newtheorem{lem}[defi]{Lemma}
\newtheorem{coro}[defi]{Corollary}
\theoremstyle{definition}
\newtheorem{rem}[defi]{Remark}
\newtheorem{ques}[defi]{Question}

\newcommand{\add}{{\rm add}}
\newcommand{\con}{{\rm con}}
\newcommand{\gd}{{\rm gl.dim}}
\newcommand{\sd}{{\rm st.dim}}
\newcommand{\sr}{{\rm sr}}
\newcommand{\dm}{{\rm dom.dim}}
\newcommand{\cdm}{{\rm codomdim}}
\newcommand{\tdim}{{\rm dim}}
\newcommand{\E}{{\rm E}}
\newcommand{\Mor}{{\rm Morph}}
\newcommand{\End}{{\rm End}}
\newcommand{\rank}{{\rm rank}}
\newcommand{\PSL}{{\rm PSL}}
\newcommand{\GL}{{\rm GL}}
\newcommand{\ind}{{\rm ind}}
\newcommand{\rsd}{{\rm res.dim}}
\newcommand{\rd} {{\rm rd}}
\newcommand{\ol}{\overline}
\newcommand{\overpr}{$\hfill\square$}
\newcommand{\rad}{{\rm rad}}
\newcommand{\soc}{{\rm soc}}
\renewcommand{\top}{{\rm top}}
\newcommand{\pd}{{\rm pdim}}
\newcommand{\id}{{\rm idim}}
\newcommand{\fld}{{\rm fdim}}
\newcommand{\Fac}{{\rm Fac}}
\newcommand{\Gen}{{\rm Gen}}
\newcommand{\fd} {{\rm fin.dim}}
\newcommand{\Fd} {{\rm Fin.dim}}
\newcommand{\Pf}[1]{{\mathscr P}^{<\infty}(#1)}
\newcommand{\DTr}{{\rm DTr}}
\newcommand{\cpx}[1]{#1^{\bullet}}
\newcommand{\D}[1]{{\mathscr D}(#1)}
\newcommand{\Dz}[1]{{\mathscr D}^+(#1)}
\newcommand{\Df}[1]{{\mathscr D}^-(#1)}
\newcommand{\Db}[1]{{\mathscr D}^b(#1)}
\newcommand{\C}[1]{{\mathscr C}(#1)}
\newcommand{\Cz}[1]{{\mathscr C}^+(#1)}
\newcommand{\Cf}[1]{{\mathscr C}^-(#1)}
\newcommand{\Cb}[1]{{\mathscr C}^b(#1)}
\newcommand{\Dc}[1]{{\mathscr D}^c(#1)}
\newcommand{\K}[1]{{\mathscr K}(#1)}
\newcommand{\Kz}[1]{{\mathscr K}^+(#1)}
\newcommand{\Kf}[1]{{\mathscr  K}^-(#1)}
\newcommand{\Kb}[1]{{\mathscr K}^b(#1)}
\newcommand{\DF}[1]{{\mathscr D}_F(#1)}

\newcommand{\Kac}[1]{{\mathscr K}_{\rm ac}(#1)}
\newcommand{\Keac}[1]{{\mathscr K}_{\mbox{\rm e-ac}}(#1)}

\newcommand{\modcat}{\ensuremath{\mbox{{\rm -mod}}}}
\newcommand{\cmodcat}{\ensuremath{\mbox{{\rm -comod}}}}
\newcommand{\Modcat}{\ensuremath{\mbox{{\rm -Mod}}}}
\newcommand{\ires}{\ensuremath{\mbox{{\rm ires}}}}
\newcommand{\Stires}{\ensuremath{\mbox{{\rm Stires}}}}
\newcommand{\Stpres}{\ensuremath{\mbox{{\rm Stpres}}}}
\newcommand{\Spec}{{\rm Spec}}

\newcommand{\stmc}[1]{#1\mbox{{\rm -{\underline{mod}}}}}
\newcommand{\Stmc}[1]{#1\mbox{{\rm -{\underline{Mod}}}}}
\newcommand{\prj}[1]{#1\mbox{{\rm -proj}}}
\newcommand{\inj}[1]{#1\mbox{{\rm -inj}}}
\newcommand{\Prj}[1]{#1\mbox{{\rm -Proj}}}
\newcommand{\Inj}[1]{#1\mbox{{\rm -Inj}}}
\newcommand{\PI}[1]{#1\mbox{{\rm -Prinj}}}
\newcommand{\GP}[1]{#1\mbox{{\rm -GProj}}}
\newcommand{\GI}[1]{#1\mbox{{\rm -GInj}}}
\newcommand{\gp}[1]{#1\mbox{{\rm -Gproj}}}
\newcommand{\gi}[1]{#1\mbox{{\rm -Ginj}}}

\newcommand{\opp}{^{\rm op}}
\newcommand{\otimesL}{\otimes^{\rm\mathbb L}}
\newcommand{\rHom}{{\rm\mathbb R}{\rm Hom}\,}
\newcommand{\pdim}{\pd}
\newcommand{\Hom}{{\rm Hom}}
\newcommand{\Coker}{{\rm Coker}}
\newcommand{ \Ker  }{{\rm Ker}}
\newcommand{ \Cone }{{\rm Con}}
\newcommand{ \Img  }{{\rm Im}}
\newcommand{\Ext}{{\rm Ext}}
\newcommand{\StHom}{{\rm \underline{Hom}}}
\newcommand{\StEnd}{{\rm \underline{End}}}
\newcommand{\KK}{I\!\!K}
\newcommand{\gm}{{\rm _{\Gamma_M}}}
\newcommand{\gmr}{{\rm _{\Gamma_M^R}}}

\def\demo{{\bf Proof}\hskip10pt}

\def\g{\gamma} \def\d{\delta} \def\a{\alpha}
\def\s{\sigma}  \def\om{\omega}  \def\ld{\lambda}
\def\D{\Delta}
\def\si{\Sigma} \def\O{\Omega}
\def\G{\Gamma} \def\GG{{\cal G}} \def\XX{{\cal X}} \def\MM{{\cal M}} \def\lrr{\lg r\rg } \def\ogg{\overline {\GG}}
\def\og{\overline G} \def\oh{\overline H} \def\oc{\overline C}
 \def\oq{\overline Q}

 \def\oa{\overline A}
  \def\ob{\overline B}  \def\ol{\overline L} \def\om{\overline M}
\def\on{\overline N} \def\op{\overline P} \def\os{\overline S}
\def\ot{\overline T} \def\ok{\overline K} \def\ov{\overline V}
\def\od{\overline D} \def\oi{\overline I} \def\oj{\overline J}
\def\o1{\overline 1} \def\olh{\overline h}

\def\o{\overline}   \def\olr{\overline r} \def\oll{\overline \ell} \def\olt{\overline t }

\def\di{\bigm|} \def\lg{\langle} \def\rg{\rangle}

\def\vez{\varepsilon}\def\bz{\bigoplus}  \def\sz {\oplus}
\def\epa{\xrightarrow} \def\inja{\hookrightarrow}

\newcommand{\lra}{\longrightarrow}
\newcommand{\llra}{\longleftarrow}
\newcommand{\lraf}[1]{\stackrel{#1}{\lra}}
\newcommand{\llaf}[1]{\stackrel{#1}{\llra}}
\newcommand{\ra}{\rightarrow}
\newcommand{\dk}{{\rm dim_{_{k}}}}

\newcommand{\holim}{{\rm Holim}}
\newcommand{\hocolim}{{\rm Hocolim}}
\newcommand{\colim}{{\rm colim\, }}
\newcommand{\limt}{{\rm lim\, }}
\newcommand{\Add}{{\rm Add }}
\newcommand{\Prod}{{\rm Prod }}
\newcommand{\pres}{\ensuremath{\mbox{{\rm pres}}}}
\newcommand{\app}{{\rm app }}
\newcommand{\Tor}{{\rm Tor}}
\newcommand{\Cogen}{{\rm Cogen}}
\newcommand{\Tria}{{\rm Tria}}
\newcommand{\Loc}{{\rm Loc}}
\newcommand{\Coloc}{{\rm Coloc}}
\newcommand{\tria}{{\rm tria}}
\newcommand{\Con}{{\rm Con}}
\newcommand{\Thick}{{\rm Thick}}
\newcommand{\thick}{{\rm thick}}
\newcommand{\Sum}{{\rm Sum}}
\def\Mon{\hbox{\rm Mon}}
\def\Aut{\hbox{\rm Aut}}
\def\SL{\hbox{\rm SL}}
\newcommand{\PGL}{{\rm PGL}}
\def\Syl{\hbox{\rm Syl}}
\def\char{\hbox{\rm \,char\,}}
{\Large \bf
\begin{center}
A classification of regular maps with Euler characteristic $-pq$
\end{center}}

\medskip
\centerline{\textbf{Xiaogang Li} and \textbf{Yao Tian} }

\medskip

\renewcommand{\thefootnote}{\alph{footnote}}
\setcounter{footnote}{-1} \footnote{}
\renewcommand{\thefootnote}{\alph{footnote}}
\setcounter{footnote}{-1}
\footnote{2020 Mathematics Subject
Classification: Primary 05C10 Secondary 20B25.}
\renewcommand{\thefootnote}{\alph{footnote}}
\setcounter{footnote}{-1}
\footnote{Keywords: Regular map; Automorphism group; Non-orientable surface; Euler characteristic.}

\begin{abstract}
In this paper, we give a classification of regular maps with Euler characteristic $-pq$ for distinct primes $q>p\geq 5$. This together with previous classification of regular maps with Euler characteristic $-2p,-3p$ and $-p^2$ completes the classification of regular maps with Euler characteristic $-pq$ for two primes $p$ and $q$. An interesting consequence is that, for every pair of twin primes $p$ and $q$ greater than $5$, there exist three regular maps with solvable automorphism groups and Euler characteristic $-pq$, up to duality and isomorphism. 
\end{abstract}

\section{Introduction}
A {\em map} is a cellular decomposition of a closed surface.
An alternative way to describe  a map is to consider it as cellular embedding of a graph $\cal{G}$ into a closed surface $\cal{S}$, such that connected components of $\cal{S}\setminus\cal{G}$ are simply connected. An {\it automorphism of a map} $\MM$ is an automorphism of the {\it underlying graph} ${\cal G}$
which extends to a  self-homeomorphism of the surface.
These automorphisms  form a subgroup $\Aut({\MM})$ of the automorphism  group $\Aut({\cal G})$ of $\cal G$.
It is well-known that  $\Aut({\MM})$ acts semi-regularly on the set of all flags (in most cases, which are incident vertex-edge-face triples).
If the action is regular, then we call the map as well as the corresponding embedding  {\it regular}.  The map is called orientable (non-orientable) if the supporting surface is orientable (non-orientable).
\vskip 2mm
There are close connections between  regular maps and  other mathematical theories such as group theory, hyperbolic geometry and complex curves and other areas. We refer the reader to \cite{JS1, Ne1, Si} for an overview. It is these connections that make classification of regular maps is an interesting and meaningful problem. During the past forty years, there has been extensive research on regular maps, see \cite{BS,CCDKNW,CDNS, GNSS,JJ,JS} for an overview.
\vskip 2mm
Classification problems of regular maps can be divided into three different approaches, that is, classifying regular maps by their supporting surfaces, underlying graphs and automorphism groups, respectively, which we refer to that as topological, graph-theoretical and group-theoretical approaches. For the second and third approaches, remarkable progress has been made in the past twenty years, for instance, a complete classification of regular maps with complete graphs as underlying graphs and descriptions of regular maps with nilpotent groups and finite simple groups as automorphism groups, see \cite{CDNS, JJ, J1, J2, Sa} for an overview. But for the topological approach, very little has been known for a long time and only those regular maps of very small genus have been completely classified \cite{MCD}. This provides evidence for the difficulties in studying regular maps from the topological approach. The first breakthrough in this direction was D'azevedo, Nedela and \v{S}ir\'{a}\v{n}'s work \cite{ARJ}, which gave a classification of regular maps on an infinite list of surfaces of Euler characteristic a negative prime. Later, Conder, Poto\v{c}nik, \v{S}ir\'{a}\v{n} and Tucker gave a classification of regular maps with Euler characteristic $-2p, -3p$ and $-p^2$ in a series of papers \cite{CNJ,MPJ,CST}. Their methods, especially their use of almost Sylow-cyclic groups (see Section 2 for definition) shed light on research in this direction. Motivated by their work, Jicheng Ma completed the classification of regular maps with Euler characteristic $-2p^2$ \cite{Ma}, and Yao Tian and Xiaogang Li successfully classified regular maps of Euler characteristic $-p^3$ \cite{TL}. 
\vskip 2mm
In this paper, we shall focus on the topological approach and we study regular maps with Euler characteristic $-pq$. Since regular maps of Euler characteristic of $-2p,-3p$ and $-p^2$ have been classified, we only need to consider the case $q>p\geq 5$ in this paper.
\vskip 2mm
Before stating our main theorem, we first introduce several groups and their related regular maps.
\vskip 2mm

Define $$\begin{array}{ll}G_1(j,k):=&\lg a,b,c,d\mid a^2=b^2=c^2=d^2=(ab)^j=1,\\
&(cd)^k=[a,c]=[a,d]=[b,c]=[b,d]=1\rg,\end{array}$$ where $2\nmid jk$. Then $G_1(j,k)\cong \mathbb{D}_{2j}\times \mathbb{D}_{2k}$.
\vskip 2mm
For any odd prime $p$, we let $\mathbb{F}_p$ denote a finite field of $p$-elements. For an even integer $n$, we define $S(n,p)\subseteq \mathbb{F}_p$ to be the set of elements $x$ in $\mathbb{F}_p$ such that the order of the following matrix 
\[M(x):=
\left(\begin{array}{cc}
0 & 1 \\
-1 & x
\end{array}\right)\in \GL(2,\mathbb{F}_p)
\]
divides $n$ but not $\frac{n}{2}$. By the natural correspondence between $\mathbb{F}_p$ and $\{0,1,2,\cdots, p-1\}$, we may also view $S(n,p)$ as a subset of $\{0,1,2,\cdots,p-1\}$. It is straightforward to check that $2\notin S(n,p)$ and  $S(n,p)= \emptyset$ if $n\leq 4$. 

Define $$\begin{array}{ll}G_2(x,n,p):=&\lg a,b,c,d\mid a^p=b^p=c^2=d^2=(cd)^n=1,\\
&[a,b]=1, a^c=a^{-1}, b^c=a^xb, a^d=b, b^d=a\rg,\end{array}$$where $x\in S(n,p)$. Then $G_2(x,n,p)\cong (\mathbb{Z}_p\times \mathbb{Z}_p)\rtimes_{\phi_x} \mathbb{D}_{2n}$, where $\phi_x$ is determined by the defining relations of $G_2(x,n,p)$. The groups $G_2(x,n,p)$ are pairwise non-isomorphic for different $x\in S(n,p)$. This can be seen from the following observations:

$(1)$ If we view $\lg a,b\rg$ as a vector space, then the conjugation of $cd$ on $\lg a,b\rg$ is given by the matrix $M(x)$ with respect to the basis $\{a,b\}$;

$(2)$ Similar matrices have the same trace;

$(3)$ All elements of order $n$ in the group $G_2(x,n,p)$ are in a conjugate calss. 

Finally, we define $$\begin{array}{ll} G_3(u):=&\lg a,b,c,d\mid a^2=b^2=[a,b]=c^2=d^2=(cd)^u=1,\\
& a^c=b,b^c=a,a^d=a,b^d=ab\rg,\end{array}$$ where $u\equiv 3({\rm mod}~6)$. Then $G_3(u)\cong \mathbb{D}_4\rtimes_\psi \mathbb{D}_{2u}$, where $\psi$ is determined by the defining relations of $G_3(u)$.

The regular maps $\MM(G_1(j,k);bc,a,d)$, $\MM(G_2(x,n,p);d,ab(cd)^{\frac{n}{2}},c)$ ($x\in S(n,p)$) and $\MM(G_3(u);c,d,ad)$ are non-orientable and of types
$\{2j,2k\}$, $\{2p,n\}$ and $\{u,4\}$, respectively, which are denoted by $\MM_1(j,k)$, $\MM_2(x,n,p)$ and $\MM_3(u)$, respectively.

\vskip 2mm
The constructions of all regular maps $\MM$ with $\Aut(M)$ being $\PSL(2,p)$ or $\PGL(2,p)$ for some prime $p\geq 5$ were given in \cite{MPJ2}. Given a regular map $\MM(\PGL(2,f);x,y,z)$ of type $\{m,n\}$, where $m$ and $n$ are divisors of $f(f^2-1)$ satisfying the conditions in \cite[Table 2]{MPJ}, if $\gcd(d,m)=1$, $x,y\in \PGL(2,f)-\PSL(2,f)$ and $z\in \PSL(2,f)$, then it is clear that $\MM(G;\alpha x,y,z)$ is a regular map of type $\{dm,n\}$ where $G:=\mathbb{Z}_d\rtimes_\varphi \PGL(2,f)$, $\mathbb{Z}_d=\lg \alpha\rg$ and $\varphi$ is the homomorphism from $\PGL(2,f)$ onto the group generated by the automorphism which maps $\alpha \in \mathbb{Z}_d$ to $\alpha^{-1}$. For convenience, we will just indicate the types and automorphism groups of such maps in this paper.
\vskip 2mm
Now we are ready to state our main theorem.
\begin{theorem}\label{main}
Let $\MM=\MM(G;r,t,\ell)$  be a regular map with Euler characteristic $-pq$ for two primes $p$ and $q$ with $q>p\geq 5$. Then one of the followings holds
\begin{enumerate}
\item[{\rm (i)}] $\MM\cong \MM_1(j,k)$ or $\MM_1(j,k)^*$, where $2\nmid jk$ and $(j-1)(k-1)=pq+1$;
\item[{\rm (ii)}] $\MM\cong \MM_2(0,4,q)$ or $\MM_2(0,4,q)^*$, where $q-p=2$;
\item[{\rm (iii)}] $\MM$ has type $\{2p^k,n\}$ and $\MM$ is a regular covering of the regular map $\MM_2(x,n,p)$ or $\MM_2(x,n,p)^*$ with transformation group a cyclic group of order $p^{k-1}$, where $n$ is even, $x\in S(n,p)$ and $p^kn-2p^k-n=2q$; 
\item[{\rm (iv)}] $\MM\cong \MM_3(u)$ or $\MM_3(u)^*$, where $u\equiv 3({\rm mod}~6)$ and $u-4=pq$;  
\item[{\rm (v)}] $\MM$ has type $\{m,n\}$ and $G\cong \PSL(2,f)$, where either \\
 {\rm (1)} $\{p,m,n\}=\{f,\frac{f-1}{2},\frac{f+1}{2}\}$ and $mn-2m-2n=2q$ or;\\
 {\rm (2)} $\{p,m,n\}=\{f,\frac{f-1}{4},\frac{f+1}{2}\}$ and $mn-2m-2n=q$ or;\\
 {\rm (3)} $\{p,m,n\}=\{f,\frac{f-1}{2},\frac{f+1}{4}\}$ and $mn-2m-2n=q$ or;\\
 {\rm (4)} $q=f$, and $p,q$ and $\{m,n\}$ are listed in Table 3.
 \item[{\rm (vi)}] $\MM$ has type $\{dm,n\}$ and $G\cong \mathbb{Z}_d\rtimes_\varphi \PGL(2,f)$, where $f$ is a prime, $p\nmid m$, $\gcd(d,f(f^2-1))\in \{1,p\}$ and either \\
 {\rm (1)} $\{p,m,n\}=\{f,f+1,\frac{f-1}{2}\}$ or $\{f,f-1,\frac{f+1}{2}\}$ and $dmn-2dm-2n=2q$, or;\\
 {\rm (2)} $d=1, p=f, \{m,n\}=\{f-1,f+1\}$ and $p^2-1-4p=4q$;, or;\\
 {\rm (3)} $d=1, p=5, q=f=7$ and $\{m,n\}=\{6,8\}$, or;\\
  {\rm (4)} $d=1, q=f$, and $p,q$ and $\{m,n\}$ are listed in Table 4.
\end{enumerate}
\end{theorem}

\begin{rem}
If $5\leq p<q$ are twin primes, i.e., $q=p+2$, then it follows from Theorem \ref{main} that $\MM_1(q,q), \MM_2(0,4,q), \MM_3(q^2-2q+4)$ and their dual are regular maps of Euler characteristic $-pq$. It is conjectured that there are infinite pairs of twin primes. If the conjecture were true, then there may exist infinitely many regular maps with Euler characteristic $-pq$ with $p,q$ being twin primes. 
\end{rem}

After this introductory section,  we will introduce some notations, terminologies, known  results and  basic theory of regular maps in  Section 2, and we will prove Theorem \ref{main} in Section 3.

\section{Preliminaries}
In this section, we shall give  some  notations and terminologies,  a brief introduction to regular maps and some known results used in this paper.

\subsection{Notations and terminologies }

For a finite set $\Omega$, the cardinality of $\Omega$ is denoted by $|\Omega|$.
For a prime $p$ and an arbitrary integer $m$, the $p$-part and $p'$-part of $m$ mean the largest $p$-power dividing $m$ and the largest divisor of $m$ coprime to $p$, respectively, and we denote them by $m_p$ and $m_{p'}$, respectively. Given a prime $p$, let $\bar{\mathbb{F}}_p$ be an algebraic closed field of characteristic $p$ and let $\mathbb{F}_p$ be its prime subfield, for a positive integer $m$, we define $\xi_m$ as a primitive $2m_{p'}$-th root of unity in $\bar{\mathbb{F}}_p$. Then we say that a set $\{m,n\}\subseteq \mathbb{N}$ is an \emph{$p$-admissible set} if $4-(\xi_m+\xi^{-1}_m)^2-(\xi_n+\xi^{-1}_n)^2$ equals to a square of an element from $\mathbb{F}_p$. 

\vskip 2mm
Let $G$ be a finite group. For $H\leq G$, $H\unlhd G$ and $H\char G$, we mean that $H$ is a subgroup, a normal subgroup and a characteristic subgroup of $G$, respectively. For $g\in G$ and $H\leq G$, by $|g|$ and $|G:H|$, we denote the order of $g$ and the index of $H$ in $G$, respectively, and by $C_G(H)$ we mean the set of elements in $G$ that commutes with all elements of $H$. For a prime $p$, the set of Sylow $p$-subgroups of $G$ (if $p\nmid |G|$, then a Sylow $p$-subgroup of $G$ is identified as the trivial subgroup $\{1\}$) and the number of Sylow $p$-subgroups of $G$ are denoted by $\Syl_p(G)$ and $n_p(G)$, respectively. The largest normal $p$-subgroup of $G$ is denoted by $O_{p}(G)$ and the largest normal subgroup of $G$ of odd order is denoted by $O_{2'}(G)$. The centre of $G$ is denoted by $Z(G)$. The Fitting subgroup of $G$ (the largest nilpotent normal  subgroups of $G$) is denoted by $F(G)$. The commutator subgroup of $G$ is denoted by $G'$. The Frattini subgroup of $G$ which is the intersection of all maximal subgroups of $G$ is denoted by $\Phi(G)$. A generating set of $G$ is a subset $T\subseteq G$ such that $G=\lg T\rg$. The group $G$ is called \emph{$d$-generated} if the minimum size of its generating set is $d$.
\vskip 2mm

For finite groups $G$ and $A$, if there is a homomorphism $\phi: A\rightarrow \Aut(G)$, then the semiproduct of $G$ and $A$ with respect to $\phi$ is denoted by $G\rtimes_\phi A$. If $G$ is abelian, then the automorphism sending every element of $G$ to its inverse is denoted by $\sigma_G$. 
\vskip 2mm
As usual, we use $\GL(n,q)$, $\PGL(n,q)$ and $\PSL(n,q)$ to denote the general linear groups, projective general linear groups and projective special linear groups of dimension $n$ over the finite field $F_q$, respectively. By $\mathbb{Z}_m$ and $\mathbb{D}_{2n}$ ($n\geq 2$), we denote the cyclic group of order $m$ and dihedral group of order $2n$, respectively.
\vskip 2mm
We introduce an important  terminology which was defined in \cite{MPJ} as follows.  A group $G$ is called {\it almost Sylow-cyclic} if all Sylow subgroups of $G$ of odd order are cyclic and all Sylow $2$-subgroup of $G$  are either trivial or contain a cyclic subgroup of index $2$.

\subsection{Regular maps}

In this subsection, we give a description of maps and regular maps in the combinatorial way. We begin with the following definition.
\begin{defi}\label{map1}
{\rm For a given finite set $F$ and three fixed-point-free involutory permutations $r, t, \ell $ on $F$, a quadruple
$\MM=\MM(F; r, t, \ell )$ is called a {\it combinatorial map} if they satisfy the following two conditions:
(1)\ $t\ell =\ell t$; (2)\ the group $\lg r,t, \ell \rg $ acts transitively on $F.$}
\end{defi}

For a given combinatorial map $\MM=\MM(F; r,t, \ell )$,  $F$  is called the {\it flag} set,
$r, t, \ell$ are called {\it rotary, transversal and longitudinal involution,} respectively.
The group $\lg r, t, \ell \rg $ is called the {\it monodromy group} of $\MM$, denoted by $\Mon(\MM)$.
We define the {\it vertices, edges} and {\it face-boundaries} of $\MM$ to be the
orbits of the subgroups $\lg r,t\rg$, $\lg t, \ell \rg $ and $\lg r, \ell \rg $, respectively.
The incidence in $\MM $ is defined by a nontrivial set intersection.

Given two maps $\MM _1=\MM(F_1; r_1, t_1, \ell _1)$ and $\MM_2=\MM_2(F_2; r_2, t_2, \ell _2)$, a bijection $\phi$ from $F_1$ to
$F_2$ is called a {\it map isomorphism} if $\phi r_1=r_2\phi$, $\phi t_1=t_2\phi$ and $\phi \ell _1=\ell _2\phi$. In particular,
if $\MM _1=\MM _2=\MM$, then $\phi$ is called an {\it automorphism} of $\MM$.
The set $\Aut(\MM)$ of all automorphisms of $\MM$ forms a group  which is called the {\it automorphism group} of $\MM$. By the definition of map isomorphism, we have $\Aut(\MM )=C_{S_F}(\Mon (\MM))$,  the centralizer of $\Mon(\MM)$ in
$S_F$, the symmetric group on $F$. It follows from the transitivity of $\Mon(\MM )$ on $F$ that $\Aut(\MM)$ acts semi-regularly on $F$.  If the action is regular, then we call $\MM$
{\it regular map}.

Let $\MM=\MM(F; r,t, \ell )$ be a regular map and fix a flag $\alpha\in F$. Then for every $\beta\in F$, there is a unique $g_{\beta}\in \Aut(\MM)$ such that $\beta=\alpha^{g_{\beta}}$ and the mapping
$$\varphi: F\rightarrow \Aut(\MM),~\beta\mapsto g_{\beta}$$
is a bijection. Let $R(\Aut(\MM))$ be the right regular representation of $\Aut(\MM)$. Then the mapping
 $$\phi: \Aut(\MM)\rightarrow R(\Aut(\MM)),~h\mapsto R(h),~\forall~ h\in \Aut(\MM)$$ is a group isomorphism.
Since $$\varphi(\beta^h)=\varphi(\alpha^{g_{\beta}h})
=g_{\beta}h=g_{\beta}^{R(h)}=\varphi(\beta)^{\phi(h)}$$
for every $\beta\in F$ and $h\in \Aut(\MM)$, the two permutation groups $\Aut(\MM)\le{\rm Sym}(F)$ and $R(\Aut(\MM))\leq {\rm Sym}(\Aut(\MM))$ are permutation isomorphic.
In this sense, we identify $F$ with the group $G=\langle r, t, \ell\rangle$ and the rotary, transversal and longitudinal involution on $F$ with the left regular representation of $r$, $t$ and $\ell$, respectively. Thus the monodromy
(automorphism) group of $\MM$ is the left (right) regular representation of $G$. Usually, we identify $R(G)$ with $G$ and call $\MM:=\MM(G; r, t, \ell )$ an {\it algebraic map}. For any normal subgroup $N\lhd G$, it is clear that $\MM(G/N; rN, tN, \ell N)$ is also an algebraic map, which is called the quotient map of $\MM$ induced by the normal subgroup $N$. On the other hand, it is obvious that $\MM^*:=\MM(G; r, \ell, t)$ is also an algebraic map, which is called the \emph{dual} of $\MM$.
It is straightforward to check that the even-word subgroup $\lg tr, r\ell \rg $ of $G$ has index at most 2.
If the index is 2, then one may fix an orientation for $\MM $ and so $\MM $ is said to be {\it orientable}.  If the index is 1, then $\MM$ is said to be {\it non-orientable}. The {\it type} of $\MM$ is defined as $\{|rt|,|r\ell|\}$, and the \emph{Euler characteristic} $\chi_{\MM}$ (resp. \emph{genus}) of $\MM$ is referred to the Euler characteristic (resp. genus) of the supporting surface of $\MM$, that is, $\chi_{\MM}:=V+F-E$. In particular,
$$
(*)\quad\quad\chi_{\MM}=-\frac{|G|(xy-2x-2y)}{4xy}$$where $x=|rt|, y=|r\ell|$. Clearly, $\chi_{\MM}=\chi_{\MM^*}$. The relation of the Euler characteristic $\chi$ and genus $g$ of a surface $S$ is given as follows:
\[
\chi=
\left\{
\begin{aligned}
& 2-2g, && S~\text{is orientable }  \\
& 2-g, && S~\text{is non-orientable }
\end{aligned}
\right.
\]
Thus only non-orientable surfaces may have odd Euler characteristic. In particular, if $\chi_{\MM}$ is odd then $\MM$ is non-orientable.

\subsection{Some Known Results}

\begin{prop}[{\cite[Chap.1, Theorem 6.11]{SUZ}}]
\label{nc}
Let $H$ be a subgroup of a group $G$.
Then $C_G(H)$ is a normal subgroup of
$N_G(H)$ and the quotient $N_G(H)/C_G(H)$ is isomorphic
to a subgroup of $\Aut (H)$.
\end{prop}

\begin{prop}[{\cite[Chap.3, Corollary 3.28]{ISA}}]
\label{cop}
Let $A$ act via automorphisms on $G$, where $A$ and $G$ are finite groups, and assume that $\gcd(|A|,|G|)=1$. If the induced action of $A$ on the Frattini factor group $G/\Phi(G)$ is trivial, then the action of $A$ on $G$ is trivial. 
\end{prop}

\begin{prop}[{\cite[Chap.5, Theorem 5.3]{ISA}}]
\label{is}
Let $G$ be a finite group and $p$ be a prime divisor of $|Z(G)\cap G'|$.
Then a Sylow-$p$ subgroup of $G$ is non-abelian.
\end{prop}

\begin{prop}[{\cite[Lemma 3.2]{MPJ}}]
\label{alm}
Let $\MM(G;r,t,\ell)$ be a regular map. If $p$ is a prime divisor of $|G|$ coprime
to the Euler characteristic of $\MM$, then a Sylow $p$-subgroup of $G$ is cyclic (if $p$ is odd) or dihedral (if $p=2$).
\end{prop}

\begin{prop}[{\cite[Theorem 1.1]{MS}}]
\label{generation}
The alternating group $A_n$ can be generated by three involutions, two of which commute if and only if $n=5$ or $n\geq 9$.
 \end{prop}

\begin{prop}[{\cite[Theorem 1]{GW}}]
\label{dihedral}
Let $G$ be a finite group with dihedral Sylow $2$-subgroups. Then $G/O_{2'}(G)$ is isomorphic to either an almost simple group with socle $\PSL(2, q)$ where $q$ is odd,
or $A_7$, or a Sylow $2$-subgroup of $G$.
 \end{prop}

\begin{prop}[{\cite[Table 1]{MPJ}}]
\label{solvable}
Let $\MM(G;r,t,\ell)$ be a non-orientable regular map with $G$ a solvable almost Sylow-cyclic group. Then we have one of the following:
\begin{enumerate}
\item[{\rm (i)}] $\MM$ has type $\{2,n\}$ and $G\cong \mathbb{D}_{2n}$, where $n\geq 4$ is even;
\item[{\rm (ii)}] $\MM$ or $\MM^*\cong \MM_1(m,n)$, where $2\nmid mn$, $1\neq m<n$ and $\gcd(m,n)=1$;
\item[{\rm (iii)}] $\MM$ or $\MM^*\cong \MM_3(m)$, where $m\equiv 3({\rm mod}~ 6)$.
\end{enumerate}
 \end{prop}

\begin{prop}[{\cite[Theorem 5.4]{MPJ}}]
\label{non-alm}
Let $\MM=\MM(G;r,t,\ell)$ be a non-orientable regular map. If $G$ is an insolvable almost Sylow-cyclic group, then one of the followings holds
\begin{enumerate}
\item[{\rm (i)}] $\MM$ has type $\{f,f\}$ and $G\cong \PSL(2,f)$, where $f\equiv 1({\rm mod}~4)$ is a prime;
\item[{\rm (ii)}] $\MM$ has type $\{m,n\}$ and $G\cong \PSL(2,f)$, where $f$ is a prime, $m\mid \frac{f\pm 1}{2}$, either $n\mid \frac{f\pm 1}{2}$ or $n=f$, and $\{m,n\}$ is an $f$-admissible set;
\item[{\rm (iii)}] $\MM$ has type $\{dm,n\}$ and  $G\cong \mathbb{Z}_d\rtimes_\varphi \PGL(2,f)$, where $f$ is a prime, $\gcd(d,f(f^2-1))=1$, $\frac{f\pm 1}{m}$ is an odd integer, and if $d>1$ then either $n=f\equiv 3({\rm mod}~4)$, or $n\mid \frac{f\pm 1}{2}$ and $\{m,n\}$ is an $f$-admissible set.
\end{enumerate}
\end{prop}

\section{Proof of Theorem~\ref{main}}

The condition of Theorem~\ref{main} is as follows:
$$(\star)~~\MM(G; r, t, \ell)~~\mbox{is a regular map with Euler characteristic}~-pq~\mbox{for primes}~q>p\geq 5.$$

Assume condition $(\star)$. Then it follows Proposition \ref{alm} that Sylow $2$-subgroups of $G$ are dihedral. Clearly, $G$ is not a dihedral group. This can be seen from \cite[Lemma 4.2]{MPJ}.  The following three facts are obvious and will be used throughout this paper without any further reference.

$(1)$ $\MM$ is non-orientable, that is, $G=\lg rt,r\ell\rg=\lg rt,t\ell\rg=\lg r\ell,t\ell\rg$;

$(2)$ $|rt|\geq 3$ and $|r\ell|\geq 3$;

$(3)$ $\lg rt\rg \cap \lg r\ell\rg=1$ (see \cite[Lemma 4.3]{MPJ}).

\vskip 2mm
Under condition $(\star)$, there are four cases $\gcd(pq, |G|)=pq$; $\gcd(pq, |G|)=1$;  $\gcd(pq, |G|)=q$; $\gcd(pq, |G|)=p$. We record them as conditions $(\star.1), (\star.2), (\star.3)$ and $(\star.4)$, respectively. The proof of Theorem \ref{main} is divided into four parts according to these four conditions. In the rest of this paper, we will always abbreviate $\MM(G;r,t,\ell)$ by $\MM$.
\vskip 2mm
The following result characterizes $\MM$ under the condition $(\star.1)$.
\begin{theorem}\label{di}
Assume condition $(\star.1)$. Then there are finitely many possibilities for $\{x,y\}, p$ and $q$.
\end{theorem}
\demo Let $\{x,y\}$ be the type of $\MM$ and let $a$ denote the $|G|$. Note that the assumption $\gcd(pq,a)=pq$ is equivalent to $pq\mid a$. This together with formula $(*)$ implies that $k(x,y)=\frac{xy}{xy-2x-2y}$ is a positive integer. It is known from \cite[Proof of Theorem 3.4]{MPJ} that there are finitely many pairs $\{x,y\}$ such that $k(x,y)$ is a positive integer, which are listed as follows:
\begin{table}[h]
\centering
\caption{List of $\{x,y\}$ and corresponding integer value of $k(x,y)$}
$
\begin{tabular}{c|c|c|c|c|c}
                 \hline
                 $\{x,y\}$ & $k(x,y)$ & $\{x,y\}$ & $k(x,y)$ & $\{x,y\}$ & $k(x,y)$ \\
                 \hline
                 \{3,7\} & 21 & \{3,24\} & 4 & \{5,5\} & 5\\

                 \{3,8\} & 12 & \{4,5\} & 10 & \{5,20\} & 2\\

                 \{3,9\} & 9 &  \{4,6\} & 6 & \{6,6\} & 3\\

                 \{3,12\}& 6 & \{4,8\}& 4 & \{6,12\}& 2 \\

                 \{3,15\}& 5 & \{4,12\}& 3 & \{8,8\}& 2\\
                 \hline
                
               \end{tabular}$
               \end{table}

Note that $x$ and $y$ are divisors of $a$ by Lagrange's theorem. With a straightforward  calculation, the possible types $\{x,y\}$ and the corresponding values of $a$ are listed in the following Table.
\begin{table}[h]
\centering
\caption{List of $\{x,y\}$ and corresponding integer value of $a$}
$
\begin{tabular}{c|c|c|c|c|c}
                 \hline
                 $\{x,y\}$ & $a$ & $\{x,y\}$ & $a$ & $\{x,y\}$ & $a$ \\
                 \hline
                 \{3,7\} & $21pq$ & \{4,5\} & $40pq$ & \{5,5\} & $20pq$\\

                 \{3,8\} & $48pq$ & \{4,6\} & $24pq$& \{6,6\} & $12pq$\\

                 \{3,9\} & $36pq$ & \{4,8\}& $16pq$ &\{8,8\}& $8pq$\\

                 \{3,12\}& $24pq$ & \{4,12\}& $12pq$ & \{5,20\}& $40q$ \\

                 \hline
                
               \end{tabular}$
               \end{table}
\vskip 2mm
Now, from Table 2, we find that $G$ is an almost Sylow-cyclic group except for the three cases: $\{x,y\}=\{3,7\}, a=84pq$ and $7\in \{p,q\}$; $\{x,y\}=\{5,5\}, a=100q$; $\{x,y\}=\{4,5\}, a=200q$. 

Suppose that $G$ is an almost Sylow-cyclic group. Then applying Proposition \ref{solvable} or Proposition \ref{non-alm} according to $G$ is solvable or not, one of the following holds
\begin{enumerate}
\item[{\rm (i)}] $\{x,y\}=\{2,n\}$ and $G\cong \mathbb{D}_{2n}$, where $n\geq 4$ is even;
\item[{\rm (ii)}] $\MM$ or $\MM^*\cong \MM_1(m,n)$, where $2\nmid mn$, $1\neq m<n$ and $\gcd(m,n)=1$;
\item[{\rm (iii)}] $\MM$ or $\MM^*\cong \MM_3(m)$, where $m\equiv 3({\rm mod}~ 6)$;
\item[{\rm (iv)}] $x=y=f$ and $G\cong \PSL(2,f)$, where $f\equiv 1({\rm mod}~4)$ is a prime;
\item[{\rm (v)}] $\{x,y\}=\{m,n\}$ and $G\cong \PSL(2,f)$, where $f$ is a prime, $m\mid \frac{f\pm 1}{2}$, either $n\mid \frac{f\pm 1}{2}$ or $n=f$, and $\{m,n\}$ is an $f$-admissible set;
\item[{\rm (vi)}] $\{x,y\}=\{dm,n\}$ and  $G\cong \mathbb{Z}_d\rtimes_\varphi \PGL(2,f)$, where $f$ is a prime, $\gcd(d,f(f^2-1))=1$, $\frac{f\pm 1}{m}$ is an odd integer, and either $n=f\equiv 3({\rm mod}~4)$, or $n\mid \frac{f\pm 1}{2}$ and $\{m,n\}$ is an $f$-admissible set.
\end{enumerate}

Case (i) is easily eliminated as $\MM$ has Euler characteristic $1\neq -pq$ in this case. 

For case (ii), we have $(m-1)(n-1)=pq+1$ and $a=4mn$.  From the equality $(m-1)(n-1)=pq+1$, we see that $mn-m-n=pq$. This implies that $p$ (resp. $q$) divides $m$ if and only if $p$ (resp. $q$) divides $n$. Thus it follows from $pq\mid a=4mn$ that $pq\mid m$ and $pq\mid n$. In particular, $pq\leq m$ and $pq\leq n$. Then $p^2q^2-2pq+1=(pq-1)(pq-1)\leq (m-1)(n-1)=pq+1$. This implies that  $pq\leq 3$,  a contradiction to our assumption that $q>p\geq 5$. 

For case (iii), we have $m-4=pq$ and $a=8m$. Since $p$ and $q$ are odd, it follows from  $pq\mid a=8m$ that $pq\mid m$. This together with $m-4=pq$ implies that $pq\mid 4$, a contradiction to our assumption that $q>p\geq 5$. 

For case (iv), we see from the Table 2 that $x=y=f=5$ and so $a=|\PSL(2,f)|=60$. Then it follows from $a=20pq$ that $pq=3$, a contradiction to our assumption that $q>p\geq 5$. 

Now assume case (v). Then we have $f(f^2-1)=8k(x,y)pq$. Clearly, $f>5$. From Table 1, we deduce that $f\notin \{x,y\}$ and so $f\mid pq$. Since $pq\mid a$ by assumption, it follows that $G\cong \PSL(2,f)$ has a nontrivial Sylow $q$-subgroup. Thus $q$ divides at least one of $f, \frac{f-1}{2}$ and $\frac{f+1}{2}$. In particular, $q\leq f$. This together with $f\mid pq$ implies that $f=q$. Hence we get $q^2-1=8k(x,y)p$. In view of Table 2, using MAGMA and Proposition \ref{non-alm} we eliminate the six cases: $\{x,y\}=\{4,12\}$ and  $\{p,q\}=\{7,13\}$; $\{x,y\}=\{4,12\}$ and $\{p,q\}=\{5,11\}$; $\{x,y\}=\{3,7\}$ and $\{p,q\}=\{11,43\}$;  $\{x,y\}=\{3,8\}$ and 
$\{p,q\}=\{23,47\}$;  $\{x,y\}=\{3,9\}$ and $\{p,q\}=\{5,19\}$; $\{x,y\}=\{4,6\}$ and $\{p,q\}=\{11,23\}$. The possibilities for $p, q$ and $\{x,y\}$ are listed in Table 3.

\begin{table}[h]
\centering
\caption{List of $\{x,y\}, p$ and $q$}
$
\begin{tabular}{c|c|c|c|c|c}
                 \hline
                 $\{x,y\}$ & $p$ & $q$ & $\{x,y\}$ & $p$ & $q$ \\
                 \hline
                 $\{3,7\}$ & 41 & 83& $\{3,12\}$ & 11 & 23\\
                 $\{3,7\}$ & 5  & 29& $\{6,6\}$  &7  & 13 \\
                 $\{3,9\}$ & 19 & 37& $\{6,6\}$ &5 & 11\\
                 \hline
                 \end{tabular}$
               \end{table}

Next, we assume  case (vi). By Proposition \ref{non-alm}, we know that $d\mid xy$. Thus, observing Table 1 we get that $d\in \{1,3,5,7,9,15\}$. Since $3\in f(f^2-1)$ and $\gcd(d,f(f^2-1))=1$, it follows that $d\notin\{ 3,9,15\}$. In particular, $d\in \{1,5,7\}$. We show that $d=1$. Suppose to the contrary that $d\neq 1$. Then $\{x,y\}=\{5,5\}$ or $\{5,20\}$ when $d=5$, and $\{x,y\}=\{3,7\}$ when $d=7$. Consider the quotient map, say $\overline{\MM}$, of $\MM$, induced by the normal cyclic subgroup of $G$ of order $d$. Then $\overline{\MM}$ is non-orientable and $\Aut(\overline{\MM})\cong \PGL(2,f)$ is a non-cyclic group. Thus $d=5$ and $\{x,y\}=\{5,20\}$. In particular, $5\cdot f(f^2-1)=40q$. This implies that $f(f^2-1)=8q$. Clearly, we have $q=f$ and $f^2-1=8$. This yields that $q=3$, a contradiction to our assumption $q>5$. Thus $d=1$. With a similar argument as for case (v), we get $q=f$ and $q^2-1=4k(x,y)p$. Now, checking Table 2 and using a case by case analysis, we have the following possibilities for $p, q$ and $\{x,y\}$: 
\begin{table}[h]
\centering
\caption{List of $\{x,y\}, p$ and $q$}
$
\begin{tabular}{c|c|c}
                 \hline
                 $\{x,y\}$ & $p$ & $q$ \\
                 \hline
                 $\{3,8\}$ & 11 & 23 \\
                 $\{3,12\}$ & 7 & 13 \\
                 $\{4,6\}$ & 7& 13\\
                 \hline
                 \end{tabular}$
               \end{table}

Finally, we consider the three cases: $\{x,y\}=\{3,7\}, a=84pq$ and $7\in \{p,q\}$; $\{x,y\}=\{5,5\}, a=100q$; $\{x,y\}=\{4,5\}, a=200q$. 

Suppose that either $\{x,y\}=\{5,5\}, a=100q$, or $\{x,y\}=\{4,5\}, a=200q$. Then $3\nmid a$. Since $3\mid |\PSL(2,f)|=\frac{f(f^2-1)}{2}$, it follows from Proposition \ref{dihedral} that $G/O_{2'}(G)$ is a 2-group. In particular, $O_{2'}(G)$ contains all elements of $G$ of odd order. Since $G=\lg rt,r\ell\rg=\lg t\ell, r\ell\rg=\lg rt, t\ell\rg$, it follows that $G=O_{2'}(G)\lg t\ell\rg$. This yields that $|G|_2=2$, contradicting $4\mid a$. 

Suppose that $\{x,y\}=\{3,7\}, a=84pq$ and $7\in \{p,q\}$. Then either $p=5, q=7$ and $a=2^2\times 3\times 5\times 7^2$, or $q>p=7$ and $a=2^2\times 3\times 7^2\times q$. According to Proposition \ref{dihedral}, $G/O_{2'}(G)$ is isomorphic to either a 2-group or an almost simple group with socle $\PSL(2,f)$. Assume that $G/O_{2'}(G)$ is isomorphic to an almost simple group with socle $\PSL(2,f)$. Comparing order, we obtain two cases: 

$(a)$ $f=5, a=2^2\times 3\times 5\times 7^2$ and $G/O_{2'}(G)\cong \PSL(2,5)$;

$(b)$ $f=q=13$ and $G/O_{2'}(G)\cong \PSL(2,13)$. 

For case $(a)$, it is clear that $O_{2'}(G)$ contains all elements of order 7. In particular, $O_{2'}(G)$ contains one of $rt$ and $r\ell$ as $\{x,y\}=\{3,7\}$. This, together with $G=\lg rt, r\ell\rg$ implies that $G/O_{2'}(G)$ is cyclic, a contradiction. 

For case $(b)$, applying Proposition \ref{is} we see that $G$ has a decomposition $G=G'\times O_{2'}(G)$, where $O_{2'}(G)\cong \mathbb{Z}_7$ and $G'\cong \PSL(2,13)$. Then $G'$ contains all elements of order $3$. This yields that $O_{2'}(G)\cong G/G'$ is cyclic of order at most $2$, a contradiction. Hence $G/O_{2'}(G)$ is isomorphic to a Sylow 2-subgroup of $G$. This implies that $O_{2'}(G)$ contains all elements of $G$ of odd order. But then $G=\lg rt, r\ell\rg\leq O_{2'}(G)<G$, a contradiction. \qed
\vskip 2mm

\subsection{The case $\gcd(pq,|G|)=1$}

In this subsection, we will determine the map $\MM$ under the condition $(\star.2)$.

\begin{theorem}\label{co}
Assume condition $(\star.2)$. Then one of the following holds
\begin{enumerate}
\item[{\rm (i)}] $\MM\cong \MM_1(j,k)$ or $\MM_1(j,k)^*$, where $\gcd(j,k)=1$ and $(j-1)(k-1)=pq+1$;
\item[{\rm (ii)}] $\MM\cong \MM_3(u)$ or $\MM_3(u)^*$, where $u\equiv 3({\rm mod}~6)$ and $u-4=pq$.
\end{enumerate}
\end{theorem}
\demo Let $a$ denote the order of $G$ and let $\{x,y\}$ be the type of $\MM$. Assume condition $(\star)$ and suppose that $\gcd(pq,|G|)=1$. According to Proposition \ref{alm}, $G$ is an almost Sylow-cyclic group. We claim that $G$ is solvable. If the claim were true, then the theorem would follow immediately from Proposition \ref{solvable}. Now, suppose to the contrary that $G$ is an insolvable group. Then by Proposition \ref{non-alm} we obtain one of the following

\begin{enumerate}
\item[{\rm (i)}] $x=y=f$ and $G\cong \PSL(2,f)$, where $f\equiv 1({\rm mod}~4)$ is a prime;
\item[{\rm (ii)}] $\{x,y\}=\{m,n\}$ and $G\cong \PSL(2,f)$, where $f$ is a prime, $m\mid \frac{f\pm 1}{2}$, either $n\mid \frac{f\pm 1}{2}$ or $n=f$, and $\{m,n\}$ is an $f$-admissible set;
\item[{\rm (iii)}] $\{x,y\}=\{dm,n\}$ and  $G\cong \mathbb{Z}_d\rtimes_\varphi \PGL(2,f)$, where $f$ is a prime, $\gcd(d,f(f^2-1))=1$, $\frac{f\pm 1}{m}$ is an odd integer, and if $d>1$ then either $n=f\equiv 3({\rm mod}~4)$, or $n\mid \frac{f\pm 1}{2}$ and $\{m,n\}$ is an $f$-admissible set.
\end{enumerate}

Thus we can write $G=\mathbb{Z}_d\rtimes N$, where $\gcd(d,f(f^2-1))=1, \PSL(2,f)\leq N\leq \PGL(2,f)$ and $N=\PGL(2,f)$ when $d>1$. Without loss of generality, assume $d\mid x$. By formula $(*)$, we obtain the equality $4pqxy=a(xy-2x-2y)$. Since $\gcd(pq,a)=1$, it follows that $pq\mid xy-2x-2y$. Thus $\frac{xy-2x-2y}{pq}=\frac{4xy}{a}$ is an integer. In particular, we have $a\leq 4xy$. Observing the three cases (i), (ii) and (iii), we find that $f\leq 5$ when $a\leq 4xy$. Then $f=5$ and $G\cong \PSL(2,5)$ or $\mathbb{Z}_d\rtimes \PGL(2,5)$. 

Assume $G\cong \PSL(2,5)$. Then $a\mid 4xy$ implies that $\{x,y\}=\{3,5\}$. However, the equality $4pqxy=a(xy-2x-2y)$ cannot hold, a contradiction. 

Assume $G\cong \mathbb{Z}_d\rtimes \PGL(2,5)$. Then we have $G/\mathbb{Z}_d\cong \PGL(2,5)$. So as the orders of elements in $\PGL(2,5)$, we have $\frac{x}{d}, y\in \{1,2,3,4,5,6\}$. Since $a\leq 4xy$, it follows that $120=\frac{a}{d}\leq 4\cdot \frac{x}{d}\cdot y$. This implies that $6\in \{\frac{x}{d},y\}$. Then $\frac{xy-2x-2y}{pq}=\frac{4xy}{a}$ is even. In particular, $2\leq \frac{4xy}{a}$, that is, $a\leq 2xy$.  But then $120=\frac{a}{d}\leq 2\cdot \frac{x}{d}\cdot y\leq 2\cdot 6^2=72$, a contradiction. \qed

\subsection{The case $\gcd(pq,|G|)=q$ }

This subsection is devoted to the characterization of $\MM$ under the conditions $(\star.3)$. 

\vskip 2mm

To proceed on, we first prove the following lemma, which characterizes $\MM(G; r, t, \ell)$ under the conditions $(\star.3)$ and $q^2\mid |G|$.
\begin{lem}\label{q^2}
Assume condition $(\star.3)$and suppose that $q^2\mid |G|$. Then one of the following holds
\begin{enumerate}
\item[{\rm (i)}] $\MM\cong \MM_1(q,q)$ or $\MM_1(q,q)^*$, where $q-2=p$;
\item[{\rm (ii)}] $\MM\cong \MM_2(0,4,q)$ or $\MM_2(0,4,q)^*$, where $q-2=p$;
\end{enumerate}
\end{lem}
\demo Let $\{x,y\}$ be the type of $\MM$. Using duality if necessary, we may assume that $x\geq y$. Recall that $y\geq 3$. Let $a$ denote $|G|$. Since the group $\lg t,\ell\rg $ of order $4$ is a subgroup of $G$, it follows that $4\mid a$. Thus by the assumption $q^2\mid a$ we may write $a=4\alpha q^2$ for an integer $\alpha$. By formula (*), we get $4\alpha q^2\cdot (xy-2x-2y)=4pqxy$. This implies that $q\mid xy$. In particular, we have $q\mid x$ or $q\mid y$. Now, we divide our analysis into three parts according to the three cases: $q\mid x$ and $q\mid y$; $q\mid x$ but $q\nmid y$; $q\nmid x$ but $q\mid y$.
\vskip 2mm
{\it Case $(1)$: $q\mid x$ and $q\mid y$.}
\vskip 2mm
In this case, we can write $x=\beta q$ and $y=\gamma q$ for integers $\beta$ and $\gamma$. Note that $\beta \geq \gamma$ by the assumption $x\geq y$. From the equality $4\alpha q^2\cdot (xy-2x-2y)=4pqxy$, we get $\alpha =\frac{p\beta \gamma}{q\beta \gamma-2\beta-2\gamma}=\frac{p\beta \gamma}{(q-2)\beta \gamma+2(\beta \gamma-\beta-\gamma)}$. If $\beta \gamma-\beta-\gamma>0$, then $\alpha<\frac{p\beta \gamma}{(q-2)\beta \gamma}\leq 1$, a contradiction to the fact that $\alpha$ is an integer. Thus $\beta \gamma-\beta-\gamma\leq 0$ and so either $\gamma=1$ or $\beta=\gamma=2$. Suppose $\gamma=1$. Then $y=q$ and $\alpha=\frac{p\beta}{(q-2)\beta-2}\leq \frac{p}{q-4}\leq \frac{q-2}{q-4}\leq 1+\frac{2}{q-4}\leq \frac{5}{3}$. This yields that $\alpha=1$, and so $a=4q^2$. By Sylow's theorem, it is easily seen that $G$ has a normal Sylow $q$-subgroup, say $Q$. Then $Q$ contains $r\ell$ and so $G=\lg t\ell, r\ell\rg\leq Q\lg t\ell\rg<G$, a contradiction. Suppose  $(\beta,\gamma)=(2,2)$. Then $\alpha=\frac{p}{q-2}$. This implies that $q-2=p$ and $\alpha=1$, and so $a=4q^2$ and $(x,y)=(2q,2q)$. By \cite[Proposition 5.3]{ARJ}, we derive that $\MM\cong  \MM_1(q,q)$ or $\MM_1(q,q)^*$, which is case (i).
\vskip 2mm
{\it Case $(2)$: $q\mid x$ and $q\nmid y$.}
\vskip 2mm
In this case, we write $x=\beta q$. Then it follows from $4\alpha q^2\cdot (xy-2x-2y)=4pqxy$ that $\alpha=\frac{p\beta y}{q\beta y-2q\beta-2y}$. Recall that $y\geq 3$. We consider the three cases $y=3$ , $y=4$ and $y\geq 5$ separately. 
\vskip 2mm
Suppose $y=3$. Then $\alpha=\frac{p\beta y}{q\beta y-2q\beta-2y}=\frac{3p\beta }{q\beta -6}$ and it follows from $y\mid a$ that $3\mid \alpha$. From $\alpha=\frac{3p\beta}{q\beta-6}$, we get $6=\beta(q-\frac{3}{\alpha}p)\geq \beta(q-p)\geq 2\beta$. Then $\beta\leq 3$ and so $\beta\in \{1,2,3\}$. This implies that $q-\frac{3}{\alpha}p$ is always an integer and so $\alpha=3$. Therefore, we have either $\beta=1$ and $q-p=6$ or $\beta=3$ and $q-p=2$. On the other hand, since $\alpha=3$, it follows that $a=12q^2$ and so $G$ has a normal Sylow $q$-subgroup $Q$ by Sylow's theorem.  Then the quotient group $G/Q$ is of order $12$ and corresponds to a regular map of type $(3,\beta)$. This implies that $\beta=3$ and $q-p=2$. But then $G/Q$ has no normal Sylow $3$-subgroup and so is isomorphic to $A_4$. This is a contradiction to the fact that $A_4$ cannot be generated by its involutions. 
\vskip 2mm
Suppose $y=4$. Then we get from $\alpha=\frac{p\beta y}{qby-2q\beta-2y}$ that $\beta(q-\frac{2}{\alpha}p)=4$. Suppose $\alpha=1$. Then $\beta(q-2p)=4$. This implies that $\beta=4$ and $q-2p=1$. Clearly, $G$ has a normal Sylow-$p$ subgroup $Q$ in this situation. By considering the quotient group of $G$ over $Q$, we see that $G/Q$ is cyclic of order $4$. But then $G/Q$ cannot be generated by its involution, a contradiction. Thus we may assume that $\alpha \geq 2$. Then $4= \beta(q-\frac{2}{\alpha}p)\geq b(q-p)\geq 2\beta$. This implies that $\beta\in \{1,2\}$. Thus $q-\frac{2}{\alpha}p$ is an integer and so $\alpha=2$. Hence we have $a=8q^2$ and either $\beta=1$ and $q-p=4$ or $\beta=2$ and $q-p=2$. Let $O$ denote the subgroup $O_{2'}(G)$. Using Proposition \ref{dihedral}, we immediately see that $G/O\cong \mathbb{D}_8$ and so $O$ is the Sylow $q$-subgroup of $G$. In particular, $G$ is solvable and $O$ contains all elements of order $q$. Suppose $\beta=1$ and $q-p=4$.  Then it follows from $G=\lg rt, t\ell\rg=\lg r\ell, t\ell\rg$ that $G/O$ is cyclic of order $2$, a contradiction. Hence $\beta=2$ and $q-p=2$. In particular, $(x,y)=(2q,4)$. Suppose that $O$ is cyclic. Then it follows from Proposition \ref{solvable} that $G$ is dihedral and $\MM$ is of Euler characteristic $1$, a contradiction. Thus $O$ is elementary abelian. Let $S$ be a Sylow 2-subgroup of $G$. Then $S\cong \mathbb{D}_8$ and $G=O\rtimes S$. By Sylow's theorem, we may choose $r,\ell$ so that $S=\lg r, \ell\rg$, and we can write $t=uv$ with $u\in O$ and $v\in S$. As $t^2=1$, we clearly have $u^v=u^{-1}$ and $v^2=1$. Since $t\ell=\ell t$, it follows that $u\ell=\ell u$ and $v\ell=\ell v$. Then $v=c\ell$, where $c$ is the unique central involution in $S$ and so $u^c=u^{-1}$. This implies that the conjugation of $S$ on $O$ is faithful as all nontrivial normal subgroups of $S$ contain $c$. Since a faithful $2$-dimensional representation over $\mathbb{F}_q$ for $\mathbb{D}_8$ is clearly irreducible and $\mathbb{D}_8$ has only one $2$-dimensional irreducible representation, the group $G$ is isomorphic to $G_2(0,4,q)$. Hence $\MM\cong \MM_2(0,4,q)$ or $\MM_2(0,4,q)^*$, which is case (ii). 
\vskip 2mm
Suppose $y\geq 5$. Then $\alpha=\frac{p\beta y}{q\beta y-2q\beta-2y}\leq \frac{p\beta y}{q\beta y-\frac{2}{5}q\beta y-\frac{2}{7}q\beta y}=\frac{35p}{11q}<\frac{35}{11}$. Thus $\alpha\in \{1,2,3\}$. Since $y\mid a=4\alpha q^2$ and $q\nmid y$, we have $y\mid 4\alpha$. Thus either $\alpha=2, y=8$ and $a=8q^2$, or $\alpha=3, y\in \{6,12\}$ and $a=12q^2$. For the case $\alpha=2, y=8, a=8q^2$, we deduce from $4\alpha q^2\cdot (xy-2x-2y)=4pqxy$ that $2p\beta=3q\beta-8$. This yields that $2p\beta=3q\beta-8=q\beta+2q\beta-8\geq 7+2(p\beta+2)-8=2p\beta+3$, a contradiction. For the case $\alpha=3, y\in \{6,12\}$ and $a=12q^2$, a contradiction can be obtained in a similar way. 

\vskip 2mm
{\it Case $(3)$: $q\nmid x$ and $q\mid y$.} 
\vskip 2mm
Since $q\mid y$, it follows that $y\geq q\geq 7$. By a similar argument as in Case (2), we deduce that $\alpha \leq 3$. Then $4\alpha \in \{4,8,12\}$. On the other hand, it follows from $q\nmid x$ and $x\mid a$ that $x\mid 4\alpha$. This, together with the assumption $x\geq y$, we have $x=4\alpha=12, y=q\in \{7, 11\}$. Substituting $x, \alpha, q$ and $y$ into the equation  $4\alpha q^2\cdot (xy-2x-2y)=4pqxy$, we get a contradiction. \qed
\vskip 2mm
Now, we are in a position to characterizes $\MM(G; r, t, \ell)$ under the conditions $(\star.3)$.
\begin{theorem}\label{q}
Assume condition $(\star.3)$. Then one of the following holds
\begin{enumerate}
\item[{\rm (i)}] $\MM\cong \MM_1(q,q)$ or $\MM_1(q,q)^*$, where $q-2=p$;
\item[{\rm (ii)}] $\MM\cong \MM_2(0,4,q)$ or $\MM_2(0,4,q)^*$, where $q-2=p$;
\item[{\rm (iii)}] $p=5, q=7$, $\MM$ has type $\{6,8\}$ and $G\cong \PGL(2,7)$. 
\end{enumerate}
\end{theorem}
\demo Suppose $q^2\mid |G|$. Then it follows from Lemma \ref{q^2} that we have either case (i) or case (ii). Thus we may assume that $q\mid |G|$ but $q^2\nmid |G|$. Using Proposition \ref{alm}, we know that $G$ is almost Sylow-cyclic. We divide our analysis into two parts according to whether $G$ is solvable or not.
\vskip 2mm
 {\it Case $(1)$: $G$ is solvable.}
\vskip 2mm
 In this case, invoking Proposition \ref{solvable}, we have one of the following cases
\begin{enumerate}
\item[{\rm (1)}] $\MM$ has type $\{2,n\}$ and $G\cong \mathbb{D}_{2n}$, where $n\geq 4$ is even;
\item[{\rm (2)}] $\MM$ or $\MM^*\cong \MM_1(m,n)$, where $2\nmid mn$, $1\neq m<n$ and $\gcd(m,n)=1$;
\item[{\rm (3)}] $\MM$ or $\MM^*\cong \MM_3(m)$, where $m\equiv 3({\rm mod}~6)$.
\end{enumerate}

The first case can be eliminated as $\MM$ has Euler characteristic $1\neq -pq$ in such case. 

For the second case, we get $|G|=4mn$ and $(m-1)(n-1)=pq+1$. Since $q\mid |G|=4mn$, it follows that $q\mid mn$ and so either $q\mid m$ or $q\mid n$. But on the other hand, we see from $(m-1)(n-1)=pq+1$ that $q\mid m$ if and only if $q\mid n$. This yields that $q^2\mid |G|$, a contradiction to our assumption. 

For the third case, we have $|G|=8m$ and $m-4=pq+1$. Then $m=pq+5$ and so $q\nmid m$ (as $q\geq 7$ by our assumption). This implies that $q\nmid 8m=|G|$, a contradiction to our assumption.
\vskip 2mm
 {\it Case $(2)$: $G$ is insolvable.}
\vskip 2mm
 In this case, using Proposition \ref{non-alm} we have one of the following
\vskip 2mm
\begin{enumerate}
\item[{\rm (a)}] $\MM$ has type $\{f,f\}$ and $G\cong \PSL(2,f)$, where $f\equiv 1({\rm mod}~4)$ is a prime;
\item[{\rm (b)}] $\MM$ has type $\{m,n\}$ and $G\cong \PSL(2,f)$, where $f$ is a prime, $m\mid \frac{f\pm 1}{2}$, either $n\mid \frac{f\pm 1}{2}$ or $n=f$, and $\{m,n\}$ is an $f$-admissible set;
\item[{\rm (c)}] $\MM$ has type $\{dm,n\}$ and  $G\cong \mathbb{Z}_d\rtimes_\varphi \PGL(2,f)$, where $f$ is a prime, $\gcd(d,f(f^2-1))=1$, $\frac{f\pm 1}{m}$ is an odd integer, and if $d>1$ then either $n=f\equiv 3({\rm mod}~4)$, or $n\mid \frac{f\pm 1}{2}$ and $\{m,n\}$ is an $f$-admissible set.
\end{enumerate}

For the first case, we have $(f-4)\cdot \frac{f-1}{2}\cdot \frac{f+1}{2}=2pq\geq 2\times 5\times 7$. This implies that $f\geq 7$. Thus, we may assume $f\geq 7$. In particular, all of $(f-4)$, $\frac{f-1}{2}$ and $\frac{f+1}{2}$ have a nontrivial prime divisors. Hence $2\in \{f-4,\frac{f-1}{2}, \frac{f+1}{2}\}$. Therefore $f\in \{3,5,6\}$, a contradiction to $f\geq 7$. 
\vskip 2mm
Assume case (b). Then we have $(mn-2m-2n)\cdot f\cdot \frac{f-1}{2}\cdot \frac{f+1}{2}=2pqmn$. This equality together with the assumption $\gcd(pq,|G|)=q$ yields that $p\mid mn-2m-2n$ and either $q=f$ or $q\mid f^2-1$. Suppose $q=f$. Then $\frac{(mn-2m-2n)}{p}\cdot \frac{f-1}{2}\cdot \frac{f+1}{2}=2mn$. We show that $f\nmid mn$. Suppose to the contrary that $f\mid mn$. Then $f\in \{m,n\}$ and so it follows from $\gcd(f,f^2-1)=1$ that $\frac{f-1}{2}\cdot \frac{f+1}{2}\leq \frac{2mn}{f}\leq 2\cdot \frac{f+1}{2}=f+1$. This yields that $f\leq 5$, a contradiction to our assumption. Also, a contradiction arises whenever both $m$ and $n$ divide $\frac{f-1}{2}$ (or $\frac{f+1}{2}$). Thus $mn$ divides $\frac{f-1}{2}\cdot \frac{f+1}{2}$ and we have two cases:

(1) $mn-2m-2n=2p$ and $\frac{f-1}{2}\cdot \frac{f+1}{2}=mn$;

(2) $mn-2m-2n=p$ and $\frac{f-1}{2}\cdot \frac{f+1}{2}=2mn$. 

Assume $mn-2m-2n=2p$ and $\frac{f-1}{2}\cdot \frac{f+1}{2}=mn$. Then we have $\{m,n\}=\{\frac{f-1}{2},\frac{f+1}{2}\}$. Substituting $m,n,q=f$ into the equation $mn-2m-2n=2p$, we get $q^2-8q-1=8p\leq 8(q-2)$. This implies that  $q<17$ and so $q\in \{7,11,13\}$. It is easily seen that the equation $q^2-8q-1=8p$ never holds. 

Assume $mn-2m-2n=p$ and $\frac{f-1}{2}\cdot \frac{f+1}{2}=2mn$. Then we have $\{m,n\}=\{\frac{f-1}{4},\frac{f+1}{2}\}$ or $\{\frac{f-1}{2},\frac{f+1}{4}\}$. Similarly, we can derive a contradiction by an argument as for the first case. Now assume that $q\mid f^2-1$. Then $q\mid \frac{f-1}{2}$ or $q\mid \frac{f+1}{2}$. In particular, we have $f\geq 2q-1$. By formula (*), we have $f\mid mn$ and so $f\mid m$ or $f\mid n$. Without loss of generality, we may assume that $f\mid m$. Thus $m=f$. Since $q^2\nmid |G|$, we derive from formula (*) that $q\nmid n$. Consequently, $4qn\mid f^2-1$, and $2=\frac{mn-2m-2n}{p}\cdot \frac{f^2-1}{4qn}$ is a product of two integers. Therefore, $(\frac{mn-2m-2n}{p},\frac{f^2-1}{4qn})=(1,2)$ or $(2,1)$. For the case $(\frac{mn-2m-2n}{p},\frac{f^2-1}{4qn})=(1,2)$, $n$ is odd and so we have $p=mn-2m-2n=(m-2)(n-2)-4\geq m-6=f-6\geq 2q-7>q+p-7\geq p$, a contradiction. For the case $(\frac{mn-2m-2n}{p},\frac{f^2-1}{4qn})=(2,1)$, $n$ is even and so $2p=mn-2m-2n\geq 2m-8=2f-8\geq 4q-10\geq 2q+4>2p$, a contradiction. 
\vskip 2mm
Finally, we assume case (c). Then $|G|=df(f^2-1)$. By the assumption $\gcd(pq,|G|)=q$, we see that $q\mid d$ or $q\mid f(f^2-1)$. Assume $q\mid d$. Then $q\nmid n$ and therefore $q^2\mid |G|$ by formula $(*)$, which is a contradiction. Thus $q\nmid d$ and so $q\mid f(f^2-1)=f(f-1)(f+1)=4f\cdot \frac{f-1}{2}\cdot \frac{f+1}{2}$. Hence $q$ divides one of $f, \frac{f-1}{2}$ and $\frac{f+1}{2}$. In particular, we have either $f=q$ or $f\geq 2q-1$. Then $|G|=df(f^2-1)=dq(q^2-1)$ when $f=q$, and $|G|=df(f^2-1)\geq d(2q-1)((2q-1)^2-1)$ when $f\neq q$. Next we show that $d=1$. Suppose to the contrary that $d>1$. Since $\gcd(d,f(f^2-1))=1$ and $6\mid f(f^2-1)$, we have $d\geq 5$. On the other hand, we have $|G|\leq 84pq$ by the Hurwitz bound on $k(x,y)=\frac{xy}{xy-2x-2y}$. Therefore, we conclude that $5q(q^2-1)\leq dq(q^2-1)\leq 84pq\leq 84q(q-2)$ when $q=f$ and $5(2q-1)((2q-1)^2-1)\leq d(2q-1)((2q-1)^2-1)\leq 84pq\leq 84q(q-2)$ when $f\neq q$. Solving the two inequalities $5q(q^2-1)\leq 84q(q-2)$ and $5(2q-1)((2q-1)^2-1)\leq 84q(q-2)$, respectively,  we get $q<17$ and $q<3$, respectively. Since we are assuming $q\geq 7$, we get $f=q<17$ and so $q\in \{7,11,13\}$ and $d\in \{5,7,11\}$. Substituting $q=f$ into the formula  $(\star)$, we get $pmn=(dmn-2m-2n)\cdot \frac{q^2-1}{4}$. Since $p\leq q-2\leq 11$, we have $11mn\geq pmn=(dmn-2m-2n)\cdot \frac{q^2-1}{4}\geq \frac{7^2-1}{4}(5mn-10m-2n)=12(mn+(4m-2)(n-2.5)-5)\geq 12(mn+10\times 0.5-5)=12mn$, a contradiction. Hence $d=1$ and so by formula (*) we have $4pqmn=f(f^2-1)(mn-2m-2n)$. Note that $q\mid f(f^2-1)$. In what follows, we consider the two cases $q\mid f$ and $q\nmid f$, separately.
\vskip 2mm
Suppose $q\neq f$. Then it follows from $q\mid f(f^2-1)$ and $q$ being odd that $q\mid \frac{f+1}{2}\cdot \frac{f-1}{2}$. This yields that $q\mid \frac{f+1}{2}$ or $q\mid \frac{f-1}{2}$. In particular, $f\geq 2q-1$. Hence $|G|=f(f^2-1)\geq 2q(2q-1)(2q-2)$. On the other hand, we know that $|G|\leq 84pq\leq 84q(q-2)$. Thus $2q(2q-1)(2q-2)\leq 84q(q-2)$. This implies that $q<11$ and so $q=7$ and $p=5$. By the inequality $|G|=f(f^2-1)\leq 84pq$ and $f\geq 2q-1$, we get $f=13$. By formula  $(\star)$, we derive that $\frac{1}{m}+\frac{1}{n}=\frac{73}{156}$. It is straightforward to check that $\frac{1}{m}+\frac{1}{n}=\frac{73}{156}$ has no integer solutions for $m$ and $n$, a contradiction. 
\vskip 2mm
Suppose $q=f$. Then $4pmn=(f^2-1)(mn-2m-2n)$ and so $p\mid (f^2-1)(mn-2m-2n)$. As $\gcd(pq,|G|)=q$, it follows that $p\nmid |G|=f(f^2-1)$ and so $p\mid (mn-2m-2n)$. Therefore, $4mn=(f^2-1)\cdot \frac{mn-2m-2n}{p}$ is a product of integers and so $f^2-1\mid 4mn$. Since $8\mid f^2-1$, it follows that $2\mid mn$. Thus $2p\mid mn-2m-2n$ and so either $f^2-1=mn$ and $mn-2m-2n=4p$, or $f^2-1=2mn$ and $mn-2m-2n=2p$. For the case $f^2-1=mn$ and $mn-2m-2n=4p$, we have $\{m,n\}=\{f-1,f+1\}$ and so $f^2-4f-1=4p$. This together with $q=f$ yields that $q^2-4q-1=4p\leq 4(q-2)$. Then $(q-1)(q-7)\leq 0$. Hence $f=q=7$, and so $p=5, \{m,n\}=\{6,8\}$,  which is case (iii) of the assertion.  For the case $f^2-1=2mn$ and $mn-2m-2n=2p$, we have $\{m,n\}=\{f-1,\frac{f+1}{2}\}$ or $\{f+1,\frac{f-1}{2}\}$. Thus either $q^2-6q-3=4p\leq 4(q-2)$ or $q^2-6q+1=4p\leq 4(q-2)$. Solving these two inequalities $q^2-6q-3\leq 4q-8$ and $q^2-6q+1\leq 4q-8$, we get $q<11$. This together with that $q\geq 7$ is a prime implies $q=7$. In particular, we have $p=5$ under the assumption  $5\leq p<q$. But then neither $q^2-6q-3=4p$ nor $q^2-6q+1=4p$ holds for $p=5$ and $q=7$, a contradiction. \qed

\subsection{The case $\gcd(pq,|G|)=p$ }

In this subsection, we are going to characterize $\MM$ under the conditions $(\star.4)$. The characterization will be done in two lemmas, one deals with the solvable case and another one deals with the insolvable case. 

\vskip 2mm
To proceed on, we need the following fact.

\begin{rem}\label{sim}
For any odd prime $f$, the general linear group $G=\GL(2,f)$ contains no simple group of the form $\PSL(2,d)$. 

To see it, suppose conversely that $G$ has such a subgroup, say $M$. Then $M$ intersects trivially with $Z=Z(\GL(2,f))$ and so $M\cong MZ/Z\leq G/Z=\PGL(2,f)$. Let $N/Z\unlhd G/Z$ be the index $2$ subgroup of $G/Z$ which is isomorphic to $\PSL(2,f)$. Then it is easy to see that $MZ/Z\leq N/Z\cong \PSL(2,f)$. Clearly, the only non-abelian simple group of $\PSL(2,f)$ is $\PSL(2,f)$ itself. Then $d=f$ and $\PSL(2,f)$ has a faithful $2$-dimension module in characteristic $f$ and so $\SL(2,f)$ has a $2$-dimension module with kernel $Z(\SL(2,f))$, which is impossible (see for instance, \cite[Chapter 1]{Alp}).
\end{rem}
Under the condition $(\star)$, using formula $(*)$ we see that any Sylow $p$-subgroup of $G$ contains a cyclic subgroup of index $p$, and thus any $p$-subgroup of $G$ is at most $2$-generated. The following lemma will be useful in our later proof.
\begin{lem}\label{act}
Let $P$ be a finite $p$-group which is at most $2$-generated, and let a non-abelian simple group $\PSL(2,f)$ act via automorphism on $P$. Then the action of $\PSL(2,f)$ on $P$ is trivial.
\end{lem}
\demo By assumption, $P$ is at most $2$-generated. Thus $P/\Phi(P)\cong \mathbb{Z}_p$ or $\mathbb{Z}_p\times \mathbb{Z}_p$. Consequently, $\Aut(P/\Phi(P))$ is isomorphic to $\mathbb{Z}_{p-1}$ or $\GL(2,p)$. Since $\PSL(2,f)$ acts via automorphism on $P$, there is a homomorphism $\phi: \PSL(2,f)\rightarrow \Aut(P)$ such that the action of $x\in \PSL(2,f)$ on $P$ is given by $\phi(x)$ on $P$. Suppose conversely that the action of $\PSL(2,f)$ on $P$ is nontrivial. Then $H:=\phi(\PSL(2,f))\cong \PSL(2,f)$. Note that, for $\sigma\in \Aut(P)$, $\sigma$ induces an automorphism, say $\sigma': P/\Phi(P)\rightarrow P/\Phi(P), x\Phi(P)\mapsto x^\sigma \Phi(P)$. This induces a homomorphism of groups $\varphi: \Aut(P)\rightarrow \Aut(P/\Phi(P)), \sigma\mapsto \sigma'$. Since $H$ is simple, we get $H\cap {\rm ker}\varphi=1$. This implies that $H\cong H{\rm ker}\varphi/{\rm ker}\varphi$ is isomorphic to a subgroup of $\Aut(P/\Phi(P))$, which forces that $\Aut(P/\Phi(P))\cong \GL(2,p)$ . But then $\GL(2,p)$ contains a subgroup which is isomorphic to $\PSL(2,f)$, a contradiction to Remark \ref{sim}. \qed
\vskip 2mm
Under the condition $(\star.4)$, the following lemma shows that $G$ can decompose as a semiproduct of a normal $p$-subgroup and an almost Sylow-cyclic subgroup.
\begin{lem}\label{semi}
Assume condition $(\star.4)$. Then $O_p(G)\in \Syl_p(O_{2'}(G))$ and $O_p(G)$ has a complement in $G$ which is an almost Sylow-cylic subgroup.
\end{lem}
\demo Let $x$ and $y$ be the order of $rt$ and $r\ell$, respectively. For convenience, we let $L$ denote $O_{2'}(G)$ and let $\pi$ be the set of prime divisors of $|L|$. Then it is clear that $O_p(G)=O_p(L)$. Set $\og:=G/O_p(L)$ and $\widehat{G}:=G/L$.

To prove that $O_p(G)\in \Syl_p(L)$, we only need to consider the case $p\in \pi$. Since $L$ is of odd order, it is solvable and so $F(L)=\prod_{s\in \pi}O_s(L)$. From $O_s(L){\rm char} L\unlhd G$, we deduce that $O_s(L)\unlhd G$ for $s\in \pi$. In particular, $C_G(O_s(L))\unlhd G$ for $s\in \pi$. By formula (*), we see that for each odd prime divisor $s$ of $|G|$ different from $p$, $|G|_s\mid x$ or $|G|_s\mid y$. This implies that either $\lg rt\rg $ or $\lg r\ell\rg$ contains a Sylow $s$-subgroup of $G$. Thus, for $s\in \pi-\{p\}$, $O_s(L)$ is contained in either $\lg rt\rg $ or $\lg r\ell\rg$ and therefore is cyclic. Hence $C_G(O_s(L))$ contains either $rt$ or $r\ell$. Since $G=\lg rt,t\ell\rg=\lg r\ell,t\ell\rg$ and $C_G(O_s(L))\unlhd G$, it follows that $G=C_G(O_s(L))\lg t\ell\rg$. In particular, $C_G(O_s(L))$ is of index at most $2$ in $G$. Hence $L\leq C_G(O_s(L))$. Suppose $O_p(L)=1$. Then $L\leq C_L(F(L))\leq F(L)\leq L$, which implies that $L=F(L)$ is cyclic. This together with $p\in \pi$ implies that $O_p(L)>1$, a contradiction. Thus $O_p(L)>1$. Note that $O_p(\ol)=1$. With a similar argument on $\og$, we derive that $\ol=F(\ol)$ is cyclic. Consequently, we see that $O_p(L)\in \Syl_p(L)$.

By Proposition \ref{dihedral}, we deduce that $\widehat{G}:=G/L$ is isomorphic to either a 2-group, $A_7$ or an almost simple group with socle $\PSL(2,f)$, where $f$ is an odd prime power. The case $\widehat{G}\cong A_7$ is ruled out by Proposition \ref{generation}, and one sees immediately that $\widehat{G}$ is isomorphic to a $2$-group or $S_4$ if $G$ is solvable, and $\widehat{G}\cong N$ for an almost simple group $N$ with socle $\PSL(2,d)$ ($d$ is an odd prime power) if $G$ is insolvable. Now, we divide our remainder of the proof into two parts according to $G$ is solvable or not.
\vskip 2mm
{\it Case $(1)$: $G$ is solvable.}
\vskip 2mm
In this case, $\widehat{G}$ is isomorphic to a $2$-group or $S_4$. Then it follows from the assumption $p\geq 5$ that $p\nmid |\widehat{G}|$, that is, $|L|_p=|G|_p$. Since $O_p(L)\in \Syl_p(L)$, it follows that $O_p(L)\in \Syl_p(G)$. Due to $O_p(L)~{\rm char}~ L\unlhd G$, we have $O_p(L)\unlhd G$. Thus it follows from the well-known Schur-Zassenhaus theorem that there exists a subgroup, say $K$, of $G$ such that $G=O_p(L)\rtimes K$. By Proposition \ref{alm}, we know that $K\cong \og$ is a solvable almost Sylow-cyclic group.  
\vskip 2mm
{\it Case $(2)$: $G$ is not solvable.}
\vskip 2mm
In this case, we have $\widehat{G}\cong N$ for an almost simple group $N$ with socle $\PSL(2,d)$, where $d$ is an odd prime power. 
\vskip 2mm
We show that $\widehat{G}$ is almost Sylow-cyclic. Suppose to the contrary that $\widehat{G}$ is not almost Sylow-cyclic. Then $\widehat{G}$ has a non-cyclic Sylow $s$-subgroup for some odd prime $s$. This implies that $d$ cannot be a prime. In particular, $d=f^a$ for an odd prime $f$ and an integer $a\geq 2$. Thus, a Sylow $f$-subgroup of $\widehat{G}$ is non-cyclic. On the other hand, using Proposition \ref{alm} we know that any Sylow $s$-subgroup of $G$ of odd order is cyclic except $s=p$. Hence, we get $f=p$ and so a Sylow $p$-subgroup of $\widehat{G}$ is elementary abelian of rank $a$. Let $x'$ and $y'$ denote $|\widehat{rt}|$ and $|\widehat{r\ell}|$, respectively. Recall that a Sylow subgroup $P\in \Syl_p(G)$ is at most $2$-generated. Consequently, a Sylow $p$-subgroup of $\widehat{G}$ is also at most two generated. This yields that $a=2$ and so $d=p^2$. In particular, $p\mid x'y'$. Without loss of generality, we may assume that $p\mid x'$. Since $(\widehat{rt})^2\in \PGL(2,p^a)$, it follows that $x'=p$ or $2p$. Now, from formula (*), we see that any odd prime divisor of $|\widehat{G}|$ divides $x'y'$. Then $y'$ divides $2(p^2-1)$ (or $2(p^2+1)$), and $2(p^2+1)$ (or $2(p^2-1)$) is a $2$-power. This forces $p=3$, a contradiction to our assumption $p\geq 5$. Hence $d$ is a prime and so $\widehat{G}$ is almost Sylow-cyclic. 
\vskip 2mm
Now, according to \cite[Theorem 1.2]{HLZZ}, we immediately see that $G$ is a split extension of $L$ by $\widehat{G}$. Consequently, there exists a subgroup, say $T$, of $G$ such that $G=L\rtimes T$. Since $O_p(L)\in \Syl_p(L)$, we know from Proposition \ref{alm} that $\og$ is also almost Sylow-cyclic. Recall that $G$ is insolvable by our assumption. This forces $\og$ to be an insolvable almost Sylow-cyclic group. According to Proposition \ref{non-alm}, we derive that either $\og \cong \PSL(2,f)$, or $\og \cong \mathbb{Z}_d\rtimes_{\varphi} \PGL(2,f)$ for some prime $f$ and an integer $d$ with $\gcd(d,f(f^2-1))=1$. 
\vskip 2mm
Suppose that $\og\cong \PSL(2,f)$. Then $L=O_p(L)$ and so $G=O_p(L)\rtimes H$ for a subgroup $H\cong \PSL(2,f)$. Note that $O_p(L)$ is at most $2$-generated. By Lemma \ref{act} we see that $G=O_p(L)\times H$. But then $H$ contains all involutions of $G$ and so $G=\lg r,t,\ell\rg\leq H$. This forces $O_p(L)=L=1$. Hence $T=G$ and we have a desired decomposition $G=O_p(L)\rtimes T$ with $T\cong \PSL(2,f)$.  
\vskip 2mm
Suppose that $\og\cong \mathbb{Z}_d\rtimes_{\varphi} \PGL(2,f)$ for some prime $f$ and an integer $d$ with $\gcd(d,f(f^2-1))=1$. Then $\widehat{G}\cong \PGL(2,f)$ and $L=O_p(L)\rtimes U$ for a subgroup $U$ of $L$ with $U\cong \mathbb{Z}_d$. Clearly, $p\nmid d$. If $O_p(L)\in \Syl_p(G)$, then it would follow from the Schur-Zassenhaus theorem that $G$ has a desired decomposition. Thus we may assume that $O_p(L)\notin \Syl_p(G)$. Then $p\mid f(f^2-1)$. To complete the proof, it suffices to show that $UT$ is a subgroup of $G$. Now, let $R\in \Syl_p(\lg rt\rg)$ and $Q\in \Syl_p(\lg r\ell\rg)$. Without loss of generality, we assume $|R|\geq |Q|$. Then observing the formula $(*)$ we deduce that $|G:R|_p\leq p$. Since $O_p(L)\notin \Syl_p(G)$, we have $O_p(L)=R$ whenever $R\leq O_p(L)$. Suppose that $O_p(L)=R$. Then $rt\in C_G(O_p(L))\unlhd G$. This, together with the basic fact $G=\lg rt,t\ell\rg$ implies that $G/C_G(O_p(L))$ is of order at most $2$. Thus $L\leq C_G(O_p(L))$ and so $L=O_p(L)\times U$ is cyclic. Thus $U~{\rm char}~L$ is normal in $G$. Hence $UT$ is a subgroup of $G$. Now, we may assume that $R\nleq O_p(L)$. Then it follows from $O_p(L)\in \Syl_p(L)$ that $p\mid x'y'$, and so any odd prime divisor of $|\widehat{G}|$ divides $x'y'$. In particular, $f\mid x'y'$. This implies that $f\in \{x',y'\}$ and either $f-1$ or $f+1$ is a $2$-power. Then either $p=f=x'$ and $y'\mid f\pm 1$, or $y'=f$ and $x'\mid f\pm 1$. On the other hand, we derive from formula $(*)$ that $|G|_{p'}\leq (4xy)_{p'}$. Since $x'=|\lg rt\rg:\lg rt\rg \cap L|$, $y'=|\lg r\ell\rg:\lg r\ell\rg \cap L|$ and $\lg rt \rg \cap \lg r\ell\rg=1$, it follows that $|\widehat{G}|_{p'}\leq (4x'y')_{p'}$. Thus, in the previous two cases, we have $p=f=5$. Considering the conjugation of $U$ on $O_p(L)$ and using Proposition \ref{nc} and Proposition \ref{cop}, we know that $U/C_U(O_p(L))$ is isomorphic to  a subgroup of either $\mathbb{Z}_{p-1}$ or $\GL(2,p)$. As $\gcd(d,f(f^2-1))=1$, we conclude that $U/C_U(O_p(L))$ is trivial, that is, $[U,O_p(L)]=1$. This implies that $L=O_p(L)\times U$. Hence $U~{\rm char}~L$ and so $U\unlhd G$. Therefore, we have $UT$ is a subgroup. \qed

\vskip 2mm
Assume condition $(\star.4)$, and suppose that $G$ is solvable. The following lemma gives a characterization of $\MM$. 
\vskip 2mm
\begin{lem}\label{so}
Assume condition $(\star.4)$ and suppose that $G$ is solvable. Then one of the following holds
\begin{enumerate}
\item[{\rm (i)}] $\MM\cong \MM_1(p^ej,p^fk)$ or $\MM_1(p^ej,p^fk)^*$, where $e$ and $f$ are positive integers with $1\in \{e,f\}$, $j, k$ are odd and $p^{e+f}jk-p^ej-p^fk=pq$;
\item[{\rm (ii)}] $\MM$ has type $\{2p^e,u\}$ and $\MM$ is a regular covering of $\MM_2(x,u,p)$ or $\MM_2(x,u,p)^*$ with transformation group a cyclic group of order $p^{e-1}$, where $e\geq 1$, $u$ is even, $x\in S(u,p)$ and $p^eu-2p^e-u=2q$.
\end{enumerate}
\end{lem}
\demo Let $x$ and $y$ be the orders of $rt$ and $r\ell$, respectively. According to Lemma \ref{semi}, we know that there exists an almost Sylow-cyclic subgroup $K$ such that $G=O_p(G)\rtimes K$. For convenience, we let $L$ and $P$ denote $O_{2'}(G)$ and $O_p(G)$, respectively, and let $Q\in Syl_p(\lg rt\rg)$ and $R\in Syl_p(\lg r\ell\rg)$. By proposition \ref{alm}, $G$ has a dihedral Sylow $2$-subgroup. Thus, using Proposition \ref{dihedral}, we conclude that $G/L$ is isomorphic to either a 2-group, $A_7$ or an almost simple group with socle $\PSL(2,f)$ where $f$ is an odd prime power. Since $G$ is solvable, we derive that $G/L$ is isomorphic to $S_4$ or a $2$-group. Thus $p\nmid |G/L|$, that is, $|G|_p=|L|_p$. By Lemma \ref{semi}, we have $P\in \Syl_p(G)$.  Since $P\unlhd G$, it follows that $Q\leq P$ and $R\leq P$. From formula $(*)$, we find that either $\frac{|G|_p}{x_p}=p$ or $\frac{|G|_p}{y_p}=p$, that is, either $|P:R|\leq p$ or $|P:Q|\leq p$. Since $\lg rt\rg \cap \lg r\ell\rg=1$, it follows that $Q\cap R=1$. Thus, if both $Q$ and $R$ are nontrivial, then we have $P=QR$ and $p\in \{|R|,|Q|\}$ by comparing the cardinalities of $P$ and $QR$. Suppose that $P$ is cyclic. Then $G$ is a solvable almost Sylow-cyclic group. Thus, using Proposition \ref{solvable} and an easy argument on the Euler characteristic of $\MM$, we find that $p\nmid |G|$, a contradiction to our assumption. Thus, $P$ is non-cyclic, and in particular $P$ is $2$-generated. Now, we divide the remainder of the proof into two parts according to whether $G/L\cong S_4$ or $G/L$ is a $2$-group. 
\vskip 2mm
{\it Case $(1)$: $G/L\cong S_4$.}
\vskip 2mm
In this case, we have $K\cong G_3(m)$, where $m\equiv 3({\rm mod}~6)$. Consider the conjugation of $K$ on $P$, which induces an action of $K$ on $P/\Phi(P)\cong \mathbb{Z}_p\times \mathbb{Z}_p$. Hence, we obtain a homomorphism of groups $\varphi: K\rightarrow \GL(2,p)$. Let $D$ be the normal subgroup of $K$ of order $4$. We show that $D\nleq {\rm ker}~\varphi$. Suppose for contradiction that $D\leq {\rm ker}~\varphi$. Then it follows from Proposition \ref{cop} that $[D,P]=1$ and so $D\unlhd G$. In particular, $|G/D|_2=2$. Consequently, $G/D$ is dihedral, and therefore $P\cong PD/D\leq G/D$ is cyclic, a contradiction. Note that any normal subgroup of $K$ either contains $D$ or intersects trivially with $D$. This, together with $D\nleq {\rm ker}~\varphi$ implies that $D\cong \varphi(D)$ is elementary abelian of order $4$. By basic linear algebra, $\varphi(D)$ is conjugate to the order $4$ group generated by the diagonal matrices ${\rm diag}(1,-1)$ and ${\rm diag}(-1,1)$ in $\GL(2,p)$. Therefore, $\varphi(D)\cap Z(\GL(2,p))$ is nontrivial. In particular, $\varphi(D)$ contains a nonidentity central element of $\varphi(K)$, which implies $Z(\varphi(K))>1$. By the first isomorphism theorem, we have $\varphi(K)\cong K/{\rm ker}~\varphi\cong G_3(m')$ for some divisor $m'$ of $m$ with $3\mid m'$. Thus $G_3(m')$ has nontrivial center, a contradiction.
\vskip 2mm
{\it Case $(2)$: $G/L$ is isomorphic to a $2$-group.}
\vskip 2mm
In this case, we have $KL/L\cong K/K\cap L$, which is a $2$-group. Note that $K\cap L\unlhd K$ has odd order. This implies that either $K\cong \mathbb{D}_{2u}$ for some even integer $u$, or $K\cong G_1(j,k)$, where $j$ and $k$ are coprime odd integers.

Suppose that $K\cong \mathbb{D}_{2u}$ for some even integer $u$. Then it follows from $P\in \Syl_p(G)$ and $G=P\rtimes K$ that $p\nmid u$. Consider the quotient map, say $\overline{\MM}$ of $\MM$, induced by the normal subgroup $P$. Then $\overline{\MM}$ has type $\{2,u\}$. Thus the type of $\MM$ is of the form $\{2p^e,up^f\}$, where $e$ and $f$ are non-negative integers. Using duality if necessary, we may assume that $x=2p^e$ and $y=up^f$. Clearly, $e\geq 1$, as $G$ is not a dihedral group. By the Schur-Zassenhaus theorem, all complements to $P$ are conjugate. Thus we may assume that $K\leq \lg r,\ell\rg$. Now we consider two cases: $f=0$ and $f\geq 1$.

Assume $f=0$. Then $\Phi(P)=\Phi(Q)$ is cyclic of order $p^{e-1}$. Set $\widetilde{G}=G/\Phi(P)$. Let $\widetilde{\MM}$ be the quotient map of $\MM$, induced by the normal subgroup $\Phi(P)$. Then $\widetilde{\MM}$ has type $\{2p,u\}$ and $\Aut(\widetilde{\MM})\cong \widetilde{G}$. Since $P$ is $2$-generated, it follows that $\widetilde{P}=P/\Phi(P)\in \Syl_p(G/\Phi(P))$ is isomorphic to $\mathbb{Z}_p\times \mathbb{Z}_p$. Clearly, $(\widetilde{rt})^2$ is a $p$-element which lies in $\widetilde{P}$. If  $\widetilde{\ell}$ normalized $\lg (\widetilde{rt})^2\rg$, then  $\lg (\widetilde{rt})^2\rg\unlhd \widetilde{G}$  and therefore the quotient map of $\widetilde{\MM}$ induced by $\lg (\widetilde{rt})^2\rg$ should have type $\{2,u\}$. Since $p\nmid u$, it follows that $P\leq \Phi(P)\lg (rt)^2\rg$. This implies that $P=\lg (rt)^2\rg$ is cyclic, a contradiction.  Thus, $\lg (\widetilde{ rt})^2\rg^{\widetilde{\ell}} \neq \lg(\widetilde{rt})^2\rg$ and so it follows from $\widetilde{P}\cong \mathbb{Z}_p\times \mathbb{Z}_p$ that $\widetilde{P}=\lg(\widetilde{rt})^2,
((\widetilde{rt})^2)^{\widetilde{\ell}} \rg$. This implies that $\widetilde{G}=\lg(\widetilde{rt})^2,(\widetilde{\ell rt\ell})^2\rg \rtimes \lg \widetilde{r}, \widetilde{\ell}\rg$. Considering the presentation of $\widetilde{G}$ using these four elements, we see that there exists $x\in S(u,p)$ such that $\widetilde{G}\cong G_2(x,u,p)$. Hence $\widetilde{\MM}\cong \MM_2(x,u,p)$ or its dual. This is Case (ii).

Assume that $f\geq 1$. Then both $Q$ and $R$ are nontrivial and so $P=QR$ and $1\in \{e,f\}$. In particular, $\Phi(P)=\Phi(Q)$ when $f=1$ and $\Phi(P)=\Phi(R)$ when $e=1$. Since $\lg rt\rg \cap \lg r\ell\rg=1$, it follows that $G=\lg rt\rg\lg r\ell\rg$ by comparing the cardinalities of both sides. By Ito's theorem, the commutator subgroup $G'$ of $G$ is abelian. Note that $(rt)^2=[r,t]$ and $(r\ell)^2=[r,\ell]$. Thus $(rt)^2,(r\ell)^2\in G'$ and so they commute. Then $[Q,R]=1$ as $Q\leq \lg (rt)^2\rg$ and $R\leq \lg (r\ell)^2\rg$. Set $T=\lg (rt)^2,(r\ell)^2 \rg$. Then $P\leq T$ and so $T=P(T\cap K)$. Note that $K\leq \lg r,\ell\rg$. Thus $T\cap K$ contains $\lg (r\ell)^2\rg \cap K$ and so $T\cap K\trianglelefteq K$. Hence $T\trianglelefteq G$, and so $|G:T|\leq 4$ as $G=\lg rt,r\ell\rg$. In particular, $J:=\lg rt,(r\ell)^2\rg \trianglelefteq G$ and $F:=\lg (rt)^2,r\ell \rg \trianglelefteq G$. Observe that $(rt)^2, (r\ell)^2\notin Z(G)$ as $(rt)^r=(rt)^{-1}, (r\ell)^r=(r\ell)^{-1}$. In particular, $[rt,(r\ell)^2]\neq 1$ and $[(rt)^2,r\ell]\neq 1$. This implies that $Z(J)=\lg (r\ell)^2\rg$ and $Z(F)=\lg (rt)^2\rg$. Then we derive from that $(rt)^t=(rt)^{-1}$ and $(r\ell)^\ell=(r\ell)^{-1}$ that $t\notin J<G$ and $\ell\notin F<G$. Suppose $t\ell\in J$ (resp. $t\ell\in F$). Then $G=\lg rt,t\ell\rg\leq J<G$ (resp. $G=\lg t\ell,r\ell\rg\leq F<G$) which is a contradiction. Thus $\ell\in J$ (resp. $t\in F$). Therefore, we deduce that $G=\lg t,(rt)^2\rg\times \lg \ell,(rt)^2\rg$, which is isomorphic to $G_1(p^e,\frac{up^f}{2})$. Hence $\MM$ is isomorphic to $\MM_1(p^e,\frac{up^f}{2})$ or its dual. This is Case (i).

Suppose that $K\cong G_1(j,k)$. We can write $K=H\rtimes D$, where $H$ is the cyclic subgroup of order $jk$ and $D$ is a Sylow $2$-subgroup of $K$. Since the type of the quotient map of $\MM$ induced by $P$ consists of two even integers, it follows that both $x$ and $y$ are even. Consider the conjugation action of $K$ on $P$, which induces an action of $K$ on $P/\Phi(P)\cong \mathbb{Z}_p\times \mathbb{Z}_p$. Hence, we obtain a homomorphism of groups $\varphi: K\rightarrow \GL(2,p)$. By proposition \ref{cop}, we have $[{\rm ker}~\varphi, P]=1$ and therefore ${\rm ker}~\varphi\unlhd G$. Suppose that $d:=|{\rm ker}~\varphi|$ is even. Then $G/{\rm ker}~\varphi$ is either cyclic or dihedral, depending on whether $4\mid d$ or not. This, together with $P\cong (P{\rm ker}~\varphi)/({\rm ker}~\varphi)\leq G/({\rm ker}~\varphi)$ implies that $P$ is cyclic, which is a contradiction. Thus $d$ is odd, and so $D\cap {\rm ker}~\varphi=1$. With a similar argument as in case (1), we deduce that $\varphi(K)$ has nontrivial center. This implies that $H\leq {\rm ker}~\varphi$. In particular, it follows from Proposition \ref{cop} that $[H,P]=1$ and hence $L=P\times H$. Then $H\unlhd G$. Note that $P\cong PH/H\leq G/H$ is non-cyclic. Thus $G/H$ can not be dihedral. In particular, $2$ does not belong to the type of the quotient map of $\MM$ induced by $H$. This yields that $p$ divides both $x$ and $y$. Consequently, both $Q$ and $R$ are nontrivial and so $P=QR$. Set $\og=G/P$. Then $\og=\lg \overline{rt}\rg\lg \overline{r\ell}\rg$. This together with $P=QR\subseteq \lg rt\rg\lg r\ell\rg$ implies that $G=\lg rt\rg\lg r\ell\rg$. With a similar argument as for the case $K\cong \mathbb{D}_{2u}$, we obtain Case (i). \qed
\vskip 2mm
Assume condition $(\star.4)$, and suppose that $G$ is an insolvable group. The following lemma gives a characterization of $\MM$. 
\begin{lem}\label{inso}
Assume condition $(\star.4)$, and suppose that $G$ is an insolvable group. Then one of the following holds
\begin{enumerate}
\item[{\rm (i)}] $\MM$ has type $\{km,n\}$ and $G\cong \mathbb{Z}_k\rtimes_\varphi \PGL(2,f)$, where $f$ is a prime, $p\nmid m$, $\gcd(k,f(f^2-1))\in \{1,p\}$ and either \\
 {\rm (1)} $\{p,m,n\}=\{f,f+1,\frac{f-1}{2}\}$ or $\{f,f-1,\frac{f+1}{2}\}$ and $kmn-2km-2n=2q$, or;\\
 {\rm (2)} $k=1, p=f, \{m,n\}=\{f-1,f+1\}$ and $p^2-4p-1=4q$.
\item[{\rm (ii)}] $\MM$ has type $\{m,n\}$ and $G\cong \PSL(2,f)$, where $f$ is a prime, and either \\
{\rm (1)} $\{p,m,n\}=\{f,\frac{f-1}{2},\frac{f+1}{2}\}$ and $mn-2m-2n=2q$ or;\\
 {\rm (2)} $\{p,m,n\}=\{f,\frac{f-1}{4},\frac{f+1}{2}\}$ and $mn-2m-2n=q$ or;\\
 {\rm (3)} $\{p,m,n\}=\{f,\frac{f-1}{2},\frac{f+1}{4}\}$ and $mn-2m-2n=q$.
\end{enumerate}
\end{lem}
\demo Let $x$ and $y$ be the order of $rt$ and $r\ell$, respectively, and let $Q\in Syl_p(\lg rt\rg)$ and $R\in Syl_p(\lg r\ell\rg)$. By Lemma \ref{semi}, we know that there exists an almost Sylow-cyclic subgroup $K$ of $G$ such that $G=O_p(G)\rtimes K$. Since $G$ is an insolvable group, it follows that $K$ is also an insolvable group, too. Let $\overline{\MM}$ be the quotient map of $\MM$ induced by $O_p(G)$. Then $\Aut(\overline{\MM})=G/O_p(G)\cong K$ is an insolvable almost Sylow-cyclic group. According to Proposition \ref{non-alm}, we have one of the following:
\begin{enumerate}
\item[{\rm (a)}] $\overline{\MM}$ has type $\{f,f\}$ and $G\cong \PSL(2,f)$, where $f\equiv 1({\rm mod}~4)$ is a prime;
\item[{\rm (b)}] $\overline{\MM}$ has type $\{m,n\}$ and $G\cong \PSL(2,f)$, where $f$ is a prime, $m\mid \frac{f\pm 1}{2}$, either $n\mid \frac{f\pm 1}{2}$ or $n=f$, and $\{m,n\}$ is an $f$-admissible set;
\item[{\rm (c)}] $\overline{\MM}$ has type $\{dm,n\}$ and  $G\cong \mathbb{Z}_d\rtimes_\varphi \PGL(2,f)$, where $f$ is a prime, $\gcd(d,f(f^2-1))=1$, $\frac{f\pm 1}{m}$ is an odd integer, and if $d>1$ then either $n=f\equiv 3({\rm mod}~4)$, or $n\mid \frac{f\pm 1}{2}$ and $\{m,n\}$ is an $f$-admissible set.
\end{enumerate}

Suppose either case (a) or case (b) holds. Note that $O_p(G)$ is at most $2$-generated. Then, by Lemma \ref{act}, $G=O_p(G)\times K$. Since $p$ is odd, it follows that $r,t,\ell\in K$. Thus $G=\lg r,t,\ell\rg\leq K$. This implies $O_p(G)=1$ and $G=K$. By direct computation, we can eliminate case (a). For case (b), by formula $(*)$ we get $f(f^2-1)(mn-2m-2n)=8pqmn$. In particular, $q\mid f(f^2-1)(mn-2m-2n)$. By the assumption $\gcd(pq,|G|)=p$, we see that $q\nmid |G|=\frac{f(f^2-1)}{2}$ and therefore $q\mid mn-2m-2n$. Thus, we obtain $f(f^2-1)\frac{mn-2m-2n}{q}=8pmn$.

Assume that $p\neq f$. Then $f\mid mn$ and it follows from $\gcd(pq,|G|)=p$ and $|G|=\frac{f(f^2-1)}{2}$ that $p\mid f^2-1$. This implies that $f\in \{m,n\}$, and $p\mid \frac{f-1}{2}$ or $p\mid \frac{f+1}{2}$. Without loss of generality, we may assume that $f=m$. Then $$(1)~~\frac{f-1}{2}\cdot \frac{f+1}{2}\cdot \frac{mn-2m-2n}{q}=2pn.$$ Note that $\gcd(\frac{f-1}{2},\frac{f+1}{2})=1$. Thus $\frac{p+1}{2}\mid 2n$ when $p\mid \frac{f-1}{2}$ and $\frac{p-1}{2}\mid 2n$ when $p\mid \frac{f+1}{2}$. Observe that if $\frac{p\pm 1}{2}$ divides $2n$ then $n$ divides $\frac{p\pm 1}{2}$. Comparing the both sides of  equation (1), we get one of the following cases \\
$(1)$ $\{p,n\}=\{\frac{f-1}{2},\frac{f+1}{2}\}$ and $mn-2m-2n=2q$;\\
$(2)$ $\{p,n\}=\{\frac{f-1}{4},\frac{f+1}{2}\}$ and $mn-2m-2n=q$;\\
$(3)$ $\{p,n\}=\{\frac{f-1}{4},\frac{f+1}{2}\}$ and $mn-2m-2n=q$.

Assume that $p=f$. Then $$(2)~~\frac{f-1}{2}\cdot \frac{f+1}{2}\cdot \frac{mn-2m-2n}{q}=2mn.$$ Clearly,  $f\nmid mn$, otherwise $f\in \{m,n\}$ and so both $\frac{f-1}{2}$ and $\frac{f+1}{2}$ divide $h\in \{m,n\}\setminus \{f\}$. This yields that $h=1$, which is impossible. Thus we derive that $m\mid \frac{p\pm 1}{2}$ and $n\mid \frac{p\mp 1}{2}$. Comparing the both sides of  equation (1), we get one of the following cases \\
$(4)$ $\{m,n\}=\{\frac{f-1}{2},\frac{f+1}{2}\}$ and $mn-2m-2n=2q$;\\
$(5)$ $\{m,n\}=\{\frac{f-1}{4},\frac{f+1}{2}\}$ and $mn-2m-2n=q$;\\
$(6)$ $\{m,n\}=\{\frac{f-1}{4},\frac{f+1}{2}\}$ and $mn-2m-2n=q$.

To summarise, we get case (ii) in the statement. 

Now, suppose that case (c) holds. Suppose $O_p(G)=1$. Then $G=K$. By a similar argument as for case (a) and case (b), we obtain case (i) for $k=d$. In the remaining proof, we consider $O_p(G)>1$. 

First, we show $G$ has a desired decomposition as in case (i). Let $N\unlhd K$ be the normal subgroup that is isomorphic to $\PSL(2,f)$, and write $K=Z\rtimes M$, where $Z\unlhd K$ is the normal subgroup of order $d$ and $M\cong \PGL(2,f)$. Then, applying Lemma \ref{act} to the conjugation of $N$ on $O_p(G)$, we see that $[O_p(G),N]=1$. Thus $N\unlhd G$. Set $\og=G/N$. Then it is clear that $|\og|_2=1$ . This yields that $\lg \olt,\oll\rg=\lg \olt\rg$ or $\lg \oll\rg$. In particular, $\og=\lg \olr,\olt,\oll\rg$ can be generated by two involutions and so is a dihedral group. Hence $\overline{O_p(G)Z}$ is cyclic. Clearly, $O_p(G)Z\cap N=1$. This implies that  $O_p(G)Z\cong O_p(G)ZN/N=\overline{O_p(G)Z}$. Thus $O_p(G)Z$ is cyclic. Set $T=O_p(G)Z$. Then $T\cong \mathbb{Z}_k$ for $k=d|O_p(G)|$ and $G=T\rtimes M$. Clearly, $M$ contains a Sylow $2$-subgroup of $G$. By Sylow theorem, we may assume that $\lg t,\ell\rg \leq M$ and we write $r=uv$ for $u\in T$ and $v\in M$. Then $G=\lg r,t,\ell\rg=\lg u,v,t,\ell\rg =\lg x\rg \rtimes \lg v,t,\ell\rg$. Comparing it with the structure of $G$, we get $T=\lg u\rg$ and $M=\lg v,t,\ell\rg$. Since $r^2=1$, $u\in T\unlhd G$ and $T\cap M=1$, it follows that $u^v=u^{-1}$ and $v^2=1$. Note that, $N\leq M$ and $[N,T]=1$. Since $u^v=u^{-1}$, it follows that $v\in M-N$ and so any element in $M-N$ inverts $T$. Hence $G\cong \mathbb{Z}_k\rtimes_\varphi \PGL(2,f)$.

We show that $\lg t,\ell\rg\nleq N$. Suppose for contradiction that $\lg t,\ell\rg\leq N$. Then both $t$ and $\ell$ commute with $u$. Thus the conjugation by the element $u^{\frac{k+1}{2}}$ maps $r=uv,t$ and $\ell$ to $v,t$ and $\ell$, respectively. This implies that $G=\lg r,t,\ell\rg \cong \lg y,t,\ell\rg\leq M$, which is clearly a contradiction by comparing their orders. Hence $\lg t,\ell\rg\nleqslant N$. Then at least one of $t$ and $\ell$ lies outside $N$. Using duality if necessary, we may assume that $t\notin N$.
Note that, for $w\in M$, $|xw|=|w|$ or ${\rm lcm}(|x|,|w|)$ depending on $w\in M-N$ or $w\in N$. Then $yt\in N$ and so $|uvt|= {\rm lcm}(|u|,|vt|)$. Recall that $1<O_p(G)\leq T=\lg u\rg$. This implies that $p\mid |uvt|$. In particular, we have $Q>1$.

We claim that $Q\leq O_p(G)$. Suppose for contradiction that $Q\nleq O_p(G)$. Note that the orders of $\overline{rt}$ and $\overline{r\ell}$ equal $m:=|vt|$ and $n:=|v\ell|$, respectively, and $O_p(G)\in \Syl_p(T)$. Set $\widehat{G}:=G/T\cong M\cong \PGL(2,f)$. Then we have $\widehat{Q}\neq 1$ and so $p\mid |\widehat{rt}|$.  From formula $(*)$, we see that $|G|_{p'}\leq (4xy)_{p'}$. Note that $m=\frac{x}{|\lg rt\rg \cap T}$ and $n=\frac{x}{|\lg r\ell\rg \cap T}$. Then it follows from $\lg rt\rg \cap \lg r\ell\rg=1$ that $|\widehat{G}|_{p'}\leq (4mn)_{p'}$. By this inequality, considering the two cases $p=f$ and $p\neq f$ separately, we find that $p=f=m=5$ and $n=6$. In particular, we see that $v\ell\in M-N$. Thus $x=|uvt|={\rm lcm}(|u|,|vt|)=k$ and $y=|uv\ell|=|v\ell|=6$. Substituting $p,f,x$ and $y$ into the formula $(*)$, we find that $q$ is even, a contradiction to the fact that $q$ is an odd prime. 

Finally, we show that $v\ell\notin N$. Suppose to the contrary that $v\ell\in N$. Since $Q\leq O_p(G)$, it follows that $Q\leq O_p(G)\leq \lg x\rg $. Note that $Q=\lg (uvt)^j\rg=\lg u^j(vt)^j\rg$, where $j={\rm lcm}(d, m_{p'})=dm_{p'}$ as $\gcd(d,m)\mid \gcd(d,f(f^2-1)=1$. Thus it follows from $Q\leq \lg x\rg$ that $m_{p'}=m$. In particular, $p\nmid m$ and so $\gcd(k,|vt|)=1$ and $x=km$. By a similar argument as in the case $Q\nleq O_p(G)$, we see that $p\nmid n$. In particular, $\gcd(k,n)=1$. Then $[v\ell,u]=1$ and it follows that $\lg uv\ell\rg=\lg u\rg\times \lg v\ell\rg$. This implies that $1<O_p(G)\leq \lg rt\rg \cap \lg r\ell\rg=1$, a contradiction. Hence $v\ell\in M\setminus N$ and so $y=|r\ell|=|uv\ell|=|v\ell|=n$. The remaining analysis is similar to that for the cases (a) and (b). \qed
\vskip 2mm
Summarizing the results of Lemma \ref{so} and Lemma \ref{inso}, we get the following theorem.

\begin{theorem}\label{p}
Assume condition $(\star.4)$. Then one of the following holds
\begin{enumerate}
\item[{\rm (i)}] $\MM\cong \MM_1(p^ej,p^fk)$ or $\MM_1(p^ej,p^fk)^*$, where $e$ and $f$ are positive integers with $1\in \{e,f\}$, $j, k$ are odd and $p^{e+f}jk-p^ej-p^fk=pq$;
\item[{\rm (ii)}] $\MM$ has type $\{2p^e,u\}$ and $\MM$ is a regular covering of $\MM_2(x,u,p)$ or $\MM_2(x,u,p)^*$ with transformation group a cyclic group of order $p^{e-1}$, where $e\geq 1$, $u$ is even and $p^eu-2p^e-u=2q$.
\item[{\rm (iii)}] $\MM$ has type $\{km,n\}$ and $G\cong \mathbb{Z}_k\rtimes_\varphi \PGL(2,f)$, where $f$ is a prime, $p\nmid m$, $\gcd(k,f(f^2-1))\in \{1,p\}$ and either \\
 {\rm (1)} $\{p,m,n\}=\{f,f+1,\frac{f-1}{2}\}$ or $\{f,f-1,\frac{f+1}{2}\}$ and $kmn-2km-2n=2q$, or;\\
 {\rm (2)} $k=1, p=f, \{m,n\}=\{f-1,f+1\}$ and $p^2-1-4p=4q$.
\item[{\rm (iv)}] $\MM$ has type $\{m,n\}$ and $G\cong \PSL(2,f)$, where $f$ is a prime, and either \\
{\rm (1)} $\{p,m,n\}=\{f,\frac{f-1}{2},\frac{f+1}{2}\}$ and $mn-2m-2n=2q$ or;\\
 {\rm (2)} $\{p,m,n\}=\{f,\frac{f-1}{4},\frac{f+1}{2}\}$ and $mn-2m-2n=q$ or;\\
 {\rm (3)} $\{p,m,n\}=\{f,\frac{f-1}{2},\frac{f+1}{4}\}$ and $mn-2m-2n=q$.
\end{enumerate}
\end{theorem}

{\bf Proof of Theorem \ref{main}} Assume condition $(\star)$. Then one of the four conditions $(\star.1), (\star.2), (\star.3)$ and $(\star.4)$ holds. Summarizing the four Theorems \ref{di}, \ref{co}, \ref{q} and \ref{p}, we get the conclusion of Theorem \ref{main}. \qed

\vskip 2mm

{\bf Acknowledgement:} This work is partially supported by the National Natural Science Foundation of China (12471332, 12350710787).

{\footnotesize
\smallskip
Xiaogang Li,
Shenzhen International Center for Mathematics, Southern University of Science and Technology, 518055
Shenzhen, Guangdong, P. R. China;

{\tt Email: 2200501002@cnu.edu.cn}

\smallskip
Yao Tian,

School of Science, Yanshan University, 066012
Qinhuangdao, Hebei, P. R. China

{\tt Email: tianyao202108@163.com}
}

\end{document}